\documentclass[pdftex,preprint]{imsart}

\usepackage{amsthm,amsmath,amssymb,natbib}
\usepackage{graphicx}

\RequirePackage[colorlinks,citecolor=blue,urlcolor=blue]{hyperref}

\startlocaldefs

\newcommand{\vt}{\vartheta}
\newcommand{\vp}{\varphi}

\newcommand{\bP}{\mathbf{P}}
\newcommand{\bp}{\mathbf{p}}
\newcommand{\bY}{\mathbf{Y}}
\newcommand{\E}{\mathbb{E}}
\newcommand{\FDR}{\mbox{FDR}}
\newcommand{\FDP}{\mbox{FDP}}
\newcommand{\FDX}{\mbox{FDX}}
\newcommand{\dis}{\stackrel{d}{=}}

\newtheorem{theorem}{Theorem}[section]
\newtheorem{lemma}{Lemma}[section]
\newtheorem{corollary}{Corollary}[section]

\endlocaldefs

\begin{document}

\begin{frontmatter}

\title{False Discovery Rate Control under Archimedean Copula}
\runtitle{Copula-based FDR Control}


\author{\fnms{Taras} \snm{Bodnar}\ead[label=e1]{bodnar@math.hu-berlin.de}}
\address{Department of Mathematics\\ Humboldt-University Berlin\\ Unter den Linden 6\\ D-10099 Berlin\\ Germany\\ \printead{e1}}
\affiliation{Department of Mathematics, Humboldt-University Berlin}
\and
\author{\fnms{Thorsten} \snm{Dickhaus}\corref{}\ead[label=e2]{dickhaus@math.hu-berlin.de}}
\address{Department of Mathematics\\ Humboldt-University Berlin\\ Unter den Linden 6\\ D-10099 Berlin\\ Germany\\ \printead{e2}}
\affiliation{Department of Mathematics, Humboldt-University Berlin}

\runauthor{Bodnar and Dickhaus}

\begin{abstract}
We are considered with the false discovery rate (FDR) of the linear step-up test $\vp^{LSU}$ considered by \cite{benhoch1995}. It is well known that $\vp^{LSU}$ controls the FDR at level $m_0 q / m$ if the joint distribution of $p$-values is multivariate totally positive of order $2$. In this, $m$ denotes the total number of hypotheses, $m_0$ the number of true null hypotheses, and $q$ the nominal FDR level. Under the assumption of an Archimedean $p$-value copula with completely monotone generator, we derive a sharper upper bound for the FDR of $\vp^{LSU}$ as well as a non-trivial lower bound. Application of the sharper upper bound to parametric subclasses of Archimedean $p$-value copulae allows us to increase the power of $\vp^{LSU}$ by pre-estimating the copula parameter and adjusting $q$. Based on the lower bound, a sufficient condition is obtained under which the FDR of $\vp^{LSU}$ is exactly equal to $m_0 q / m$, as in the case of stochastically independent $p$-values. Finally, we deal with high-dimensional multiple test problems with exchangeable test statistics by drawing a connection between infinite sequences of exchangeable $p$-values and Archimedean copulae with completely monotone generators. Our theoretical results are applied to important copula families, including Clayton copulae and Gumbel copulae.
\end{abstract}


\tableofcontents
\end{frontmatter}

\section{Introduction} \label{sec1}
Control of the false discovery rate (FDR) has become a standard type I error criterion in large-scale multiple hypotheses testing. When the number $m$ of hypotheses to be tested simultaneously is of order $10^3 - 10^6$, as it is prevalent in many modern applications from the life sciences like genetic association analyses, gene expression studies, functional magnetic resonance imaging, or brain-computer interfacing, it is typically infeasible to model or to estimate the full joint distribution of the data. Hence, one is interested in generic procedures that control the FDR under no or only qualitative assumptions regarding this joint distribution.
The still by far most popular multiple test for FDR control, the linear step-up test $\vp^{LSU}$ (say) considered in the seminal work by \cite{benhoch1995}, operates on marginal $p$-values $p_1, \hdots, p_m$. As shown by \cite{benyek:2001} and \cite{Sarkar2002}, $\vp^{LSU}$ is generically FDR-controlling over the class of models that lead to positive dependency among the random $p$-values $P_1, \hdots, P_m$ in the sense of positive regression dependency on subsets (PRDS), including $p$-value distributions which are multivariate totally positive of order $2$ (MTP$_2$). Under the PRDS assumption, the FDR of $\vp^{LSU}$ is upper-bounded by $m_0 q / m$, where $m_0$ denotes the number of true null hypotheses and $q$ the nominal FDR level.

In this work, we extend these findings by deriving a sharper upper bound for the FDR of $\vp^{LSU}$ in the case that the dependency structure among $P_1, \hdots, P_m$ can be expressed by an Archimedean copula. Our respective contributions are threefold. First, we quantify the magnitude of conservativity (non-exhaustion of the FDR level $q$) of $\vp^{LSU}$ in various copula models as a function of the copula parameter $\eta$. This allows for a gain in power in practice by pre-estimating $\eta$ and adjusting the nominal value of $q$. Second, we demonstrate by computer simulations that the proposed upper bound leads to a robust procedure in the sense that the variance of this bound over repeated Monte Carlo simulations is much smaller than the corresponding variance of the false discovery proportion (FDP) of $\vp^{LSU}$. This makes the utilization of our upper bound an attractive choice in practice, addressing the issue that the FDP is typically not well concentrated around its mean, the FDR, if $p$-values are dependent. As a by-product, we directly obtain that the FDR of $\vp^{LSU}$ is bounded by $m_0 q / m$ under the assumption of an Archimedean $p$-value copula, without explicitly relying on the MTP$_2$ property (which is fulfilled in the class of Archimedean $p$-value copulae with completely monotone generator functions, cf. \cite{Scarsini2005}). Let us point out already here that the FDR criterion is only suitable if the number $m$ of tests is large. In this case, the restriction to completely monotone generators is essentially void, because every copula generator is necessarily $m$-monotone. Third, in an asymptotic setting ($m \to \infty$), we show that the class of Archimedean $p$-value copulae with completely monotone generators includes certain models with $p$-values or test statistics, respectively, which are exchangeable under null hypotheses, $H_0$-exchangeable for short. Such $H_0$-exchangeable test statistics occur naturally in many multiple test problems, for instance in many-to-one comparisons or if test statistics are given by jointly Studentized means (cf.\ \cite{FiDi4711}).

In addition, we also derive and discuss a lower FDR bound for $\vp^{LSU}$ in terms of the generator of an Archimedean $p$-value copula. 
Application of this lower bound leads to sufficient conditions under which the FDR of $\vp^{LSU}$  is exactly equal to $m_0 q / m$, at least asymptotically as $m$ tends to infinity and $m_0/m$ converges to a fixed value. Hence, if the latter conditions are fulfilled, the FDR behaviour of $\vp^{LSU}$ is under dependency the same as in the case of jointly stochastically independent $p$-values.
 
The paper is organized as follows. In Section \ref{sec2}, we set up the necessary notation, define our class of statistical models for $P_1, \hdots, P_m$, and recall properties and results around the FDR. Our main contributions are presented in Section \ref{sec3}, dealing with FDR control of $\vp^{LSU}$ under the assumption of an Archimedean copula. Special parametric copula families are studied in Section \ref{sec4}, where we quantify the realized FDR of $\vp^{LSU}$ as a function of $\eta$. Section \ref{sec5} outlines methods for pre-estimation of $\eta$. We conclude with a discussion in Section \ref{sec6}. Lengthy proofs are deferred to Section \ref{sec_proofs}.

\section{Notation and preliminaries} \label{sec2}
All multiple test procedures considered in this work depend on the data only via (realized) marginal $p$-values $p_1, \hdots, p_m$ and their ordered values $p_{(1)}\le p_{(2)} \le \hdots \le p_{(m)}$.
Hence, it suffices to model the distribution of the random vector $\bP=(P_1, \hdots, P_m)^\top$ of $p$-values
and we consider statistical models of the form
$([0, 1]^m, \mathcal{B}([0, 1]^m), (\mathbb{P}_{\vt, \eta}: \vt \in \Theta, \eta \in \Xi))$.
In this, we assume that $\vt$ is the (main) parameter of statistical interest and we identify the null hypotheses
$H_{i} : 1 \leq i \leq m$ with non-empty subsets of $\Theta$, with corresponding alternatives $K_i = \Theta \setminus H_i$. The null hypothesis $H_i$ is called true if $\vt \in H_i$ and false otherwise. We let
$I_0 \equiv I_0(\vt) = \{1 \leq i \leq m: \vt \in H_i\}$ denote the index set of true hypotheses and $m_0 \equiv m_0(\vt) = |I_0|$ the number of true nulls. Without loss of generality, we will assume $I_0(\vt) = \{1, \hdots, m_0\}$ throughout the work. Analogously, we define $I = \{1, \hdots, m\}$, $I_1 \equiv I_1(\vt) = I \setminus I_0$ and $m_1 \equiv m_1(\vt) = |I_1| =  m - m_0$.
The intersection hypothesis $H_0 = \bigcap_{i=1}^m H_i$ will be referred to as the global (null) hypothesis.

The parameter $\eta$ is the copula parameter of the joint distribution of $\bP$, thus representing the dependency structure among $P_1, \hdots, P_m$. Its parameter space $\Xi$ may be of infinite dimension. In particular, in Section \ref{sec3} we will consider the class of all Archimedean copulas which can be indexed by the generator function $\psi$. However, we sometimes restrict our attention to parametric subclasses, for instance the class of Clayton copulae which can be indexed by a one-dimensional copula parameter $\eta \in \mathbb{R}$. In any case, we will assume that $\eta$ is a nuisance parameter in the sense that it does not depend on $\vt$ and that the marginal distribution of each $P_i$ is invariant with respect to $\eta$. Therefore, to simplify notation, we will write $\mathbb{P}_\vt$ instead of
$\mathbb{P}_{\vt, \eta}$ if marginal $p$-value distributions are concerned. Throughout the work, the $p$-values
$P_1, \hdots, P_m$ are assumed to be valid in the sense that
$$
\forall 1 \leq i \leq m: \forall \vt \in H_i: \forall t \in [0, 1]: \mathbb{P}_\vt(P_i \leq t) \leq t.
$$

A (non-randomized) multiple test operating on $p$-values is a measurable mapping $\vp = (\vp_i: 1 \leq i \leq m) : [0, 1]^m \to \{0, 1\}^m$ the components of which have the usual interpretation of a statistical test for $H_i$ versus $K_i$, $1 \leq i \leq m$. For fixed $\vp$, we let $V_m \equiv V_m(\vt) = |\{i \in I_0(\vt): \vp_i = 1\}|$ denote the (random) number of false rejections (type I errors) of $\vp$ and  $R_m \equiv R_m(\vt) = |\{i \in \{1, \hdots, m\}: \vp_i = 1\}|$ the total number of rejections.
The FDR under $(\vt, \eta)$ of $\vp$ is then defined by

\begin{equation}\label{FDR}
\text{FDR}_{\vt, \eta}(\vp) =\E_{\vt, \eta}\left[\left(\frac{V_m}{R_m \vee 1}\right)\right],
\end{equation}
and $\vp$ is said to control the FDR at level $q \in (0, 1)$ if $\sup_{\vt \in \Theta, \eta \in \Xi} \text{FDR}_{\vt, \eta}(\vp) \leq q$. The random variable $V_m / \max(R_m, 1)$ is referred to as the false discovery proportion of $\vp$, 
$\FDP_{\vt, \eta}(\vp)$ for short. Notice that, although the trueness of the null hypotheses is determined by $\vt$ alone, the FDR depends on $\vt$ and $\eta$, because the dependency structure among the $p$-values typically influences the distribution of $\vp$ when regarded as a statistic
with values in $\{0, 1\}^m$.

The linear step-up test $\vp^{LSU}$, also referred to as Benjamini-Hochberg test or \textsl{the} FDR procedure in the literature, rejects exactly hypotheses $H_{(1)}, \hdots, H_{(k)}$, where the bracketed indices correspond to the order of the $p$-values and $k=\max\{1 \leq i \leq m: p_{(i)} \le q_i\}$ for linearly increasing critical values $q_i=iq/m$. If $k$ does not exist, no hypothesis is rejected.
The sharpest characterization of FDR control of $\vp^{LSU}$ that we are aware of so far is given in the following theorem.

\begin{theorem}[\cite{AORC}] \label{das-FDR-theorem} $ $\\
Consider the following assumptions.
\begin{description}
\item[(D1)] $\forall (\vt, \eta) \in \Theta \times \Xi: \forall j\in I: \forall i\in I_0(\vt)$:  $\mathbb{P}_{\vt, \eta}(R_m \geq j | P_i \leq t)$ is non-increasing in $t\in(0, q_j]$.
\item[(D2)] $\forall \vt\in\Theta: \forall i\in I_0(\vt): P_i \sim \mbox{UNI}([0,1])$.
\item[(I1)]$\forall (\vt, \eta) \in \Theta \times \Xi$: The $p$-values $(P_i: i\in I_0(\vt))$ are independent and identically distributed (iid).
\item[(I2)] $\forall (\vt, \eta) \in \Theta \times \Xi$: The random vectors $(P_i: i\in I_0(\vt))$ and $(P_i: i\in I_1(\vt))$ are stochastically independent.
\end{description}
Then, the following two assertions hold true.
\begin{eqnarray}
\text{Under (D1),~} \forall (\vt, \eta) \in \Theta \times \Xi: \quad \FDR_{\vt, \eta}(\varphi^{LSU}) &\leq &\frac{m_0(\vt)}{m} q. \label{FDR-control-LSU}\\
\text{Under (D2)-(I2),~}
\forall (\vt, \eta) \in \Theta \times \Xi: \quad \FDR_{\vt, \eta}(\varphi^{LSU}) &= &\frac{m_0(\vt)}{m} q. \label{FDR-equal-LSU}
\end{eqnarray}
\end{theorem}
The crucial assumption (D1) is fulfilled for multivariate distributions of $\bP$ which are positively regression dependent on the subset $I_0$ (PRDS on $I_0$) in the sense of \cite{benyek:2001}. In particular, if the joint distribution of $\bP$ is MTP$_2$, then (D1) holds true.

To mention also a negative result, \cite{guo-rao-JSPI} have shown that there exists a multivariate distribution of $\bP$ such that the FDR of $\vp^{LSU}$ is equal to $m_0 q / m \sum_{j=1}^m j^{-1}$, showing that $\vp^{LSU}$ is not generically FDR-controlling over all possible joint distributions of $\bP$. The main purpose of the present work (Section \ref{sec3}) is to derive a sharper upper bound on the right-hand side of \eqref{FDR-control-LSU}, assuming that $\Xi$ is the space of completely monotone generator functions of Archimedean copulae.

\section{FDR control under Archimedean Copula} \label{sec3}

In this section, it is assumed that the joint distribution of $\bP$ is given by an Archimedean copula such that
\begin{equation}\label{P-joint_dis}
F_{\bP}(p_1,...,p_m)=\mathbb{P}_{\vt, \psi}(P_1\le p_1,...,P_m\le p_m) = \psi\left(\sum_{i=1}^m \psi^{-1}\left(F_{P_i}(p_i)\right)\right),
\end{equation}
where the function $\psi(\cdot)$ is the so-called copula generator and takes the role of $\eta$ in our general setup. In \eqref{P-joint_dis} and throughout the work, $F_\xi$ denotes the cumulative distribution function (cdf) of the variate $\xi$. 
The generator $\psi$ fully determines the type of the Archimedean copula; see, e.g. \cite{nelsen2006}. A necessary and sufficient condition under which a function $\psi: \mathbb{R}_+ \rightarrow [0,1]$ with $\psi(0)=1$ and $\lim_{x \rightarrow \infty}\psi(x)=0$ can be used as a copula generator is that $\psi(\cdot)$ is an $m$-altering function, that is, $(-1)^d \psi^{(d)}(\cdot) \ge 0$ for $d \in \{1,2,...,m\}$, cf. \cite{Scarsini2005}. Throughout the present work, we impose a slightly stronger assumption on $\psi$. Namely, we assume that $\psi$ is completely monotone, i.\ e. $(-1)^d \psi^{(d)}(\cdot) \ge 0$ for all $d \in \mathbb{N}$. If $m$ is large as it is usual in applications of the FDR criterion, the distinction between the class of completely monotone functions and the class of $m$-altering functions becomes negligible.

A very useful property of an Archimedean copula with completely monotone generator $\psi$ is the stochastic representation of $\bP$. Namely, there exists a sequence of jointly independent and identically UNI$[0,1]$-distributed random variables $Y_1, \hdots, Y_m$ such that (cf. \cite{Marshall-Olkin}, Section 5)
\begin{equation}\label{P-stoch_pres}
\bP = (P_i: 1 \leq i \leq m) \dis \left(F^{-1}_{P_i}\left(\psi\left(\log\left(Y_i^{-1/Z}\right)\right)\right): 1 \leq i \leq m \right),
\end{equation}
where the symbol $\dis$ denotes equality in distribution.
The random variable $Z$ with Laplace transform $t \mapsto \psi(t)=\E[(\exp(-tZ))]$ is independent of $Y_1, \hdots, Y_m$,
and its distribution is determined by $\psi$ only. Throughout the remainder, $\mathbb{P}$ and $\mathbb{E}$ refer to the distribution of $Z$, for ease of presentation. The stochastic representation (\ref{P-stoch_pres}) shows that the type of the Archimedean copula can equivalently be expressed in terms of the random variable $Z$. Moreover, the $p$-values $(P_i: 1 \leq i \leq m)$ are conditionally independent given $Z=z$. This second property allows us to establish the following sharper upper bound for the FDR of $\vp^{LSU}$.

\begin{theorem}[Upper FDR bound] \label{coro-sharp}
Let $Z$ be as in \eqref{P-stoch_pres} and
let $\bP^{(i)}$ consist of the $(m-1)$ remaining p-values obtained by dropping $P_i$ from $\bP$ so that $P^{(i)}_{(1)}\le P^{(i)}_{(2)}\le ...\le P^{(i)}_{(m-1)}$. The random set $D_k^{(i)}$ is then given by
\begin{equation} \label{D_k-i-P}
D_k^{(i)} = \{q_{k+1}\le P_{(k)}^{(i)}, \hdots, q_{m}\le P_{(m-1)}^{(i)} \}.
\end{equation}
For a given value $Z=z$ we define the function $\mathbf{T}: [0,1]^m \rightarrow [0,1]^m$ by
$\mathbf{T}(\bp)=(T_1(p_1),...,T_m(p_m))^T$ with $T_j(p_j)=\exp \left(-z\psi^{-1}\left(F_{P_j}(p_j)\right) \right)$
for $\bp=(p_1,...,p_m)^T \in [0,1]^m$. This function transforms, for fixed $Z=z$, realizations of $\bP$ into realizations of
$\bY = (Y_1, \hdots, Y_m)^\top$ given in \eqref{P-stoch_pres}. Let $D_{\bY;k}^{(i,z)}$ denote the image of the set $D_{k}^{(i)}$ under $\mathbf{T}$ for given $Z=z$ and let $G_k^i(z) = \mathbb{P}_{\vt, \psi} \left(D_{\bY;k}^{(i,z)}\right)$.
Then it holds
$$
\forall \vt \in \Theta: \FDR_{\vt, \psi}(\varphi^{LSU}) \leq \frac{m_0(\vt)}{m}q - b(m, \vt, \psi),
$$
where
\begin{eqnarray} \label{sharper_upper_bound}
b(m, \vt, \psi) \!\!\! &=& \!\!\! \frac{q}{m} \sum_{i=1}^{m_0} \sum_{k=1}^{m-1} \int_{z^*_k}^\infty \left(\frac{\exp \left(-z\psi^{-1}(q_{k+1})\right)}{q_{k+1}} - \frac{\exp\left(-z\psi^{-1}(q_{k})\right)}{q_k}\right) \times \nonumber\\
& &\hspace*{100pt} (G_k^i(z)-G_k^i(z^*_k)) dF_Z(z) \nonumber\\
\!\!\! & =& \!\!\! \frac{q}{m} \sum_{i=1}^{m_0} \sum_{k=1}^{m-1}
\E\left[ \left(\frac{\exp \left(-Z\psi^{-1}(q_{k+1})\right)}{q_{k+1}} - \frac{\exp\left(-Z\psi^{-1}(q_{k})\right)}{q_k}\right) \times \right. \nonumber\\
& &\hspace*{100pt} \left. (G_k^i(Z)-G_k^i(z^*_k)) \mathbf{1}_{[z^*_k,\infty)} (Z) \right] 
\end{eqnarray}
with
\begin{equation}\label{z^*_k-Sec2}
z^*_k=\frac{\log q_{k+1}-\log q_{k}}{\psi^{-1}(q_{k})-\psi^{-1}(q_{k+1})}=\frac{\log\left(1+1/k\right)}{\psi^{-1}(kq/m)-\psi^{-1}((k+1)q/m)}
\end{equation}
and $\mathbf{1}_{A}$ denoting the indicator function of the set $A$.
\end{theorem}

Noticing that $b(m, \vt, \psi)$ is always non-negative, we obtain the following result
as a straightforward corollary of Theorem \ref{coro-sharp}.

\begin{corollary}\label{thm1}
Let the copula of $\bP=(P_1,...,P_m)^\top$ be an Archimedean copula,
where $P_i$ is continuously distributed on $[0,1]$ for $1 \leq i \leq m$. Then it holds that
\begin{equation}\label{FDR_upper}
\forall \vt \in \Theta: \forall \psi \in \Xi: \FDR_{\vt, \psi}(\varphi^{LSU}) \le \frac{m_0(\vt)}{m} q,
\end{equation}
where $\Xi$ denotes the set of all completely monotone generator functions of Archimedean copulae.
\end{corollary}

The result of Corollary \ref{thm1} is in line with the findings obtained by \cite{benyek:2001} and \cite{Sarkar2002} that we have recalled in Section \ref{sec1}. Namely, \cite{Scarsini2005} pointed out that an Archimedean copula possesses the MTP$_2$ property if the copula generator $\psi$ is completely monotone and, hence, the FDR is controlled by $\varphi^{LSU}$ in this case.


From the practical point of view, it is problematic that $b(m, \vt, \psi)$ depends on the (main) parameter $\vt$ of statistical interest. In practice, one will therefore often only be able to work with $\sup_{\vt \in \Theta} \{m_0(\vt) q / m - b(m, \vt, \psi)\}$. Since $b(m, \vt, \psi) \geq 0$ for all $\vt \in \Theta$, the latter $\vt$-free upper bound will typically still yield an improvement over the "classical" upper bound. The issue of minimization of $b(m, \cdot, \psi)$ over $\vt \in \Theta$ is closely related to the challenging task of determining least favorable parameter configurations (LFCs) for the FDR.
So-called Dirac-uniform configurations are least favorable (provide upper FDR bounds) for $\vp^{LSU}$ under independence assumptions and are assumed to be generally least favorable for $\vp^{LSU}$ also in models with dependent $p$-values, at least for large values of $m$ (cf., e.\ g., \cite{FiDi4711}, \cite{Sinica2013}). \cite{Troendle2000} motivated the consideration of
Dirac-uniform configurations from the point of view of consistency of marginal tests with respect to the sample size. 
Furthermore, the expectations in \eqref{sharper_upper_bound} can in general not be calculated analytically. However, they can easily be approximated by means of computer simulations. Namely, the approximation is performed by generating random numbers which behave like independent realizations of $Z$, which completely specifies the type of the Archimedean copula, evaluating the functions $G_k^i$ at the generated values and replacing the theoretical expectation of $Z$ by the arithmetic mean of the resulting values of the integrand in \eqref{sharper_upper_bound}. Under Dirac-uniform configurations, evaluation of $G_k^i$ can efficiently be performed by means of recursive formulas for the joint cdf of the order statistics of $\bY$. We discuss these points in detail in Section \ref{sec4}.

Next, we discuss a lower bound for the FDR of $\vp^{LSU}$ under the assumption of an Archimedean copula. 

\begin{theorem}[Lower FDR bound] \label{theorem-lower-bound}
Let the copula of $\bP=(P_1,...,P_m)^T$ be an Archimedean copula with generator function $\psi$, where $P_i$ is continuously distributed on $[0,1]$ for $i=1, \hdots, m$. Then it holds that

\begin{equation}\label{FDR_lower}
\forall \vt \in \Theta: \FDR_{\vt, \psi}(\varphi^{LSU}) \geq \frac{m_0 q}{m} \gamma_{min},
\end{equation}
where
\begin{equation} \label{gamma-min1}
\gamma_{min} \equiv \gamma_{min}(\psi) = \int \min_{k \in \{1,...,m\}} \left\{\frac{\exp\left(-z\psi^{-1}\left(kq / m\right)\right)}{kq/m}\right\} dF_Z(z).
\end{equation}
\end{theorem}

From the assertion of Theorem \ref{theorem-lower-bound} we conclude that the lower bound for the FDR of $\varphi^{LSU}$ under the assumption of an Archimedean copula crucially depends on the extreme points of the function $g(\cdot | z)$, given by
\begin{equation} \label{function-g}
g(x|z)=\frac{\exp\left(-z x\right)}{\psi(x)}
\end{equation}
for $x \in \{ \psi^{-1}(q/m), \psi^{-1}(2q/m), \hdots, \psi^{-1}(q)\}$. If for all $z>0$ the minimum of $g(x|z)$ is always attained for the same index $k^*$ (say), then $\gamma_{min}=1$ and together with Theorem \ref{thm1} we get $\FDR_{\vt, \psi}(\varphi^{LSU})= m_0(\vt) q / m$. This follows directly from the identity
$$
\int \frac{\exp\left(-z \psi^{-1}\left(k^* q / m\right)\right)}{k^* q/m}  dF_{Z}(z) =
\frac{\psi \left( \psi^{-1}\left(k^* q / m \right)\right)}{k^* q/m} = 1.
$$

However, the latter holds true only in some specific cases. To obtain a more explicit constant $\gamma_{min}(\psi)$ in the general case, we notice that, due to the analytic properties of $\psi$, there exists a point $z^*$ such that $g\left(\psi^{-1}(q)|z\right)<g\left(\psi^{-1}(q/m)|z\right)$ for $z<z^*$ and $g\left(\psi^{-1}(q)|z\right)>g\left(\psi^{-1}(q/m)|z\right)$ for $z>z^*$.
The point $z^*$ is obtained as the solution of
\begin{eqnarray*}
0 &= &g\left(\psi^{-1}(q)|z\right)-g\left(\psi^{-1}(q/m)|z\right)\\
&= &\frac{\exp\left(-z \psi^{-1}\left(q\right)\right)}{q}-\frac{\exp\left(-z \psi^{-1}\left(q / m\right)\right)}{q/m}\\
&=& \frac{1}{q}\left(\exp\left(-z \psi^{-1}\left(q\right)\right)-\exp\left(\log m -z \psi^{-1}\left(\frac{q}{m}\right)\right)\right)\\
&=& \frac{\exp\left(-z \psi^{-1}\left(q\right)\right)}{q}\left(1-\exp\left(\log m + z \left(\psi^{-1}\left(q\right)-\psi^{-1}\left(\frac{q}{m}\right)\right)\right)\right),
\end{eqnarray*}
which leads to
\begin{equation}\label{z-star}
z^*=\frac{\log m}{\psi^{-1}\left(q /m \right)-\psi^{-1}\left(q\right)}\,.
\end{equation}

Next, we analyze the function $x \mapsto g(x|z)$ for given $z$. For its derivative with respect to $x$, it holds that
\begin{eqnarray*}
g^\prime(x|z)&=&-\frac{\exp\left(-z x\right)}{(\psi(x))^2}\left(z\psi(x)+\psi^\prime(x)\right).
\end{eqnarray*}
Setting this expression to zero, we get that any extreme point of $g(\cdot|z)$
 satisfies
\begin{equation}\label{x_z}
z\psi(x)+\psi^\prime(x)=0.
\end{equation}

Let $x_z$ be a solution of (\ref{x_z}). Then, the second derivative of $g(\cdot|z)$ at $x_z$ is given by
\begin{equation}
g^{\prime\prime}(x_z|z) =  -\frac{\exp\left(-z x_z\right)}{(\psi(x_z))^2}\left(z\psi^\prime(x_z)+\psi^{\prime \prime}(x_z)\right). \label{g_der2}
\end{equation}
Substituting \eqref{x_z} with $x = x_z$ in \eqref{g_der2}, we obtain
\begin{eqnarray*}
g^{\prime\prime}(x_z|z)&=& -\frac{\exp\left(-z x_z\right)}{(\psi(x_z))^2}\left(-\frac{(\psi^\prime(x_z))^2}{\psi(x_z)}+\psi^{\prime \prime}(x_z)\right)\\
&=&-\frac{\exp\left(-z x_z\right)}{(\psi(x_z))^3}\left(\psi(x_z)\psi^{\prime \prime}(x_z)-(\psi^\prime(x_z))^2\right)
\end{eqnarray*}
and application of the Cauchy-Schwarz inequality leads to
\begin{eqnarray*}
\psi(x_z)\psi^{\prime \prime}(x_z)&=& \int \exp\left(-z x_z\right) dF_Z(z)\int z^2 \exp\left(-z x_z \right) dF_Z(z)\\
&\ge&\left(\int z \exp\left(-z x_z \right) dF_Z(z)\right)^2=(\psi^\prime(x_z))^2.
\end{eqnarray*}
This proves that $g^{\prime\prime}(x_z|z)\le 0$ if $x_z$ is an extreme point of $g(x_z|z)$. Thus, any such $x_z$ is a maximum and the minimum in \eqref{gamma-min1} is attained at $\psi^{-1}(q)$ for $z \le z^*$ as well as at $\psi^{-1}(q/m)$ for $z \ge z^*$. This allows for a more explicit characterization of the lower bound.

\begin{lemma} \label{gamma-min-lem}
The quantity $\gamma_{min} \equiv \gamma_{min}(\psi)$ from \eqref{gamma-min1}
can equivalently be expressed as
\begin{eqnarray}
\gamma_{min} &= &1-\int_0^{z^*} \left(g\left(\psi^{-1}(q/m)|z\right)-g\left(\psi^{-1}(q)|z\right) \right)dF_Z(z) \label{FDR_ineq1}\\
&=& 1-\E\left[(g\left(\psi^{-1}(q/m)|Z\right)-g\left(\psi^{-1}(q)|Z\right)) \mathbf{1}_{[0,z^*]}(Z)\right], \label{FDR_ineq2}
\end{eqnarray}
where $g(\cdot|z)$ and $z^*$ are defined in \eqref{function-g} and \eqref{z-star}, respectively.
\end{lemma}

If the integral in \eqref{FDR_ineq1} cannot be calculated analytically, then it can easily be approximated via a Monte Carlo simulation by using the expression on the right-hand side of \eqref{FDR_ineq2}
and replacing the theoretical expectation by its pseudo-sample analogue.

Lemma \ref{gamma-min-lem} possesses several interesting applications. We consider the quantity $\gamma_{min}$ itself. It holds that $1 \geq \gamma_{min} \geq \underline{\gamma}_{min}$, where
\begin{equation}\label{underline_gamma}
\underline{\gamma}_{min}=1-\min \left\{\int_0^{z^*} \sup_{z\in[0,z^*]} h(z) dF_Z(z), \int_{z^*}^\infty \sup_{z\in[z^*,\infty]}(-h(z)) dF_Z(z)  \right\}
\end{equation}
with
\begin{eqnarray*} 
h(z) &= &g\left(\psi^{-1}(q/m)|z\right)-g\left(\psi^{-1}(q)|z\right)\\
&= &\frac{\exp\left(-z \psi^{-1}\left(\frac{q}{m}\right)\right)}{q/m}-\frac{\exp\left(-z \psi^{-1}\left(q\right)\right)}{q},
\end{eqnarray*}
because
\begin{equation} \label{symmetry}
\int_0^{z^*} h(z) dF_Z(z)=-\int_{z^*}^\infty h(z) dF_Z(z).
\end{equation}
However, both of the integrals in \eqref{symmetry} can be bounded by different values.
To see this, we study the behavior of the function $z \mapsto h(z)$. It holds that
\begin{eqnarray*}
h^\prime(z)&=& -\psi^{-1}\left(\frac{q}{m}\right)\frac{\exp\left(-z \psi^{-1}\left(\frac{q}{m}\right)\right)}{q/m}+\psi^{-1}\left(q\right)\frac{\exp\left(-z \psi^{-1}\left(q\right)\right)}{q}\\
&=&\psi^{-1}\left(q\right)\frac{\exp\left(-z \psi^{-1}\left(q\right)\right)}{q} \times\\
& &\hspace*{5pt} \left(1-\exp\left(\log m +\log \frac{\psi^{-1}\left(\frac{q}{m}\right)}{\psi^{-1}\left(q\right)} + z \left(\psi^{-1}\left(q\right)-\psi^{-1}\left(\frac{q}{m}\right)\right)\right)\right).
\end{eqnarray*}
Since $\psi^{-1}$ is a non-increasing function, we get that there exists a unique minimum of $h(z)$ at
\[
\underline{z}^*=\frac{\log m + \log \psi^{-1}\left(q/m\right) - \log \psi^{-1}\left(q\right)}{\psi^{-1}\left( q/m \right)-\psi^{-1}\left(q\right)}\ge z^*.
\]
Consequently, we get
\begin{eqnarray*}
\int_0^{z^*} \sup_{z\in[0,z^*]} h(z) dF_Z(z) &=&\int_0^{z^*} h(0) dF_Z(z)=h(0)F_Z(z^*)\\
&= &\frac{m-1}{q}F_Z(z^*),\\
 \int_{z^*}^\infty \sup_{z\in[z^*,\infty]}(-h(z)) dF_Z(z)&=& \int_{z^*}^\infty h(\underline{z}^*) dF_Z(z)=h(\underline{z}^*)(1-F_Z(z^*))\\
 &=&\frac{\exp\left(-\underline{z}^* \psi^{-1}\left(q\right)\right)}{q}\left(1-\frac{\psi^{-1}\left(q\right)}{\psi^{-1}\left(q/m\right)}\right) \times\\
 & & \; \; \; (1-F_Z(z^*)).
\end{eqnarray*}

\begin{corollary} \label{out-of-support}
Under the assumptions of Theorem \ref{theorem-lower-bound}, the following two assertions hold true.

\begin{enumerate}
\item[(a)] If $z^*$ from \eqref{z-star} does not lie in the support of $F_Z$, i.\ e., if $F_Z(z^*)=0$ or $F_Z(z^*)=1$, then $\gamma_{min}=\underline{\gamma}_{min}=1$
and, consequently, $\FDR_{\vt, \psi}(\varphi^{LSU}) = m_0 q / m$.
\item[(b)] Assume that $\pi_0 = \lim_{m \to \infty} m_0 / m$ exists. If $z^* = z^*(m)$
is such that $F_Z\left(z^*(m)\right) \rightarrow 0$ or $F_Z\left(z^*(m)\right) \rightarrow 1$ as $m \to \infty$, then
$$
\lim_{m \to \infty} \FDR_{\vt, \psi}(\varphi^{LSU}) = \pi_0 q.
$$
\end{enumerate}
\end{corollary}

Part (b) of Corollary \ref{out-of-support} motivates a deeper consideration of asymptotic or high-dimensional multiple tests,
i.\ e., the case of $m \rightarrow \infty$, under our general setup. This approach has already been discussed
widely in previous literature. For instance, it was called "asymptotic multiple test" by \cite{genowass:2002}. The case $m \rightarrow \infty$ was also considered by \cite{finnrot1998}, \cite{storey:2002a}, \cite{genowass2004}, \cite{FiDi4711, AORC}, \cite{JinCai2007}, \cite{SunCai2007}, and \cite{CaiJin2010}, among others.

Very interesting connections can be drawn between Archimedean $p$-value copulae and infinite sequences of $H_0$-exchangeable $p$-values. More precisely, let us assume an infinite sequence $(P_i)_{i \in \mathbb{N}}$ of $p$-values which are absolutely continuous and uniformly distributed on $[0, 1]$ under the respective null hypothesis $H_i$. Furthermore, we let $F_i$ denote the cdf. of $P_i$ under $\vt$ and assume that
$F_1(P_1), \hdots, F_m(P_m), \hdots$ are exchangeable random variables, entailing that $P_1,\hdots$, $P_m, \hdots$ themselves are exchangeable under the global hypothesis $H_0$. Sequences of $H_0$-exchangeable $p$-values have already been investigated by \cite{finnrot1998} and \cite{FiDi4711} in special settings. Moreover, the assumption of exchangeability is also pivotal in other areas of statistics, let us mention Bayesian analysis and validity of permutation tests. 
The problem of exchangeability in population genetics has been discussed by \cite{Kingman1978}.

For ease of notation, let $\tilde{P}_i= F_i(P_i)$ for $i \in \mathbb{N}$.
Because $\tilde{P}_1,\hdots,\tilde{P}_m,\hdots$ is an exchangeable sequence of random variables, it exists a random variable $Z$ with distribution function $F_Z$ such that the joint distribution of $\tilde{P}_1,\hdots,\tilde{P}_{m}$ is for any fixed $m \in \mathbb{N}$ given by 
\begin{equation}\label{th1_eq1}
F_{\tilde{P}_1,\hdots,\tilde{P}_{m}}(p_1,\hdots,p_m)=\int F_{\tilde{P}_1|Z=z}(p_1) \times
\hdots \times F_{\tilde{P}_{m}|Z=z}(p_{m}) dF_Z(z),
\end{equation}
see \cite{Olshen1974} and equation (3.1) of \cite{Kingman1978}.
Moreover, assuming that $Z \in (0, \infty)$ with probability $1$, we obtain for any $i \in \mathbb{N}$ from \cite{Marshall-Olkin}, p. 834, that
$$ 
p_i=F_{\tilde{P}_i}(p_i)=\int \exp\left(-z \psi^{-1}(p_i)\right)dF_Z(z)\,,
$$ 
where $\psi$ denotes the Laplace transform of $Z$, i.\ e., $\psi(t)=\mathbb{E}[\exp(-tZ)]$.

Theorem \ref{thm-ex} establishes a connection between the finite-dimensional marginal distributions of $H_0$-exchangeable $p$-value sequences and Archimedean copulae.

\begin{theorem} \label{thm-ex}
Assume that the elements in the infinite sequence $(P_i)_{i \in \mathbb{N}}$ are absolutely continuous and $H_0$-exchangeable. Furthermore, let the following two assumptions be valid.
\begin{itemize}
\item[(i)] The random variable $Z$ from \eqref{th1_eq1} 
takes values in $(0, \infty)$ with probability $1$.
\item[(ii)] It holds
\begin{equation} \label{equation-21}
F_{\tilde{P}_i|Z=z}(p_i)= \exp\left(-z \psi^{-1}(p_i)\right), z \in (0, \infty).
\end{equation}
\end{itemize}
Then, for any $m$,
\[
\mathbf{p} = (p_1, \hdots, p_m)^\top \mapsto \psi \left(\sum_{i=1}^{m} \psi^{-1}(p_i)\right)
\]
is a copula of $P_1,...,P_{m}$, where $\psi(t)=\mathbb{E}[\exp(-tZ)]$.
\end{theorem}

The final result of this section is an immediate consequence of
Theorem \ref{thm-ex} and Corollary \ref{thm1}.

\begin{corollary}
Under the assumptions of Theorem \ref{thm-ex}, it holds:
\begin{enumerate}
\item[a)] Any $m$-dimensional marginal distribution of the sequence 
$(P_i)_{i \in \mathbb{N}}$ possesses the MTP$_2$ property, $m \geq 2$.
\item[b)] The linear step-up test $\vp^{LSU}$, 
applied to $p_1, \hdots, p_m$, controls the FDR at level $q$.
\end{enumerate}
\end{corollary}

\section{Examples: Parametric copula families} \label{sec4}

In this section, we apply the theoretical results of Section \ref{sec3} to several parametric families of Archimedean copulas.

\subsection{Independence Copula}
The generator of the independence copula is given by
$\psi(t)=\exp(-t)$. Substituting $\psi^{-1}(x)= - \ln (x)$
in \eqref{gamma-min1}, we get
\[\gamma_{min}=  \min_{k \in \{1,...,m\}} \left\{\frac{\exp\left(\ln \left(kq / m \right)\right)}{kq/m}\right\}=1\]
and, hence,
\[
\forall \vt \in \Theta: FDR_{\vt, \psi}(\vp^{LSU}) = \frac{m_0(\vt)}{m} q
\]
under the assumption of independence. This result is in line with the previous finding
reported in \eqref{FDR-equal-LSU}.

\subsection{Clayton Copula}
\begin{figure}[ptbh]
\begin{center}%
\includegraphics[width=10.7cm]{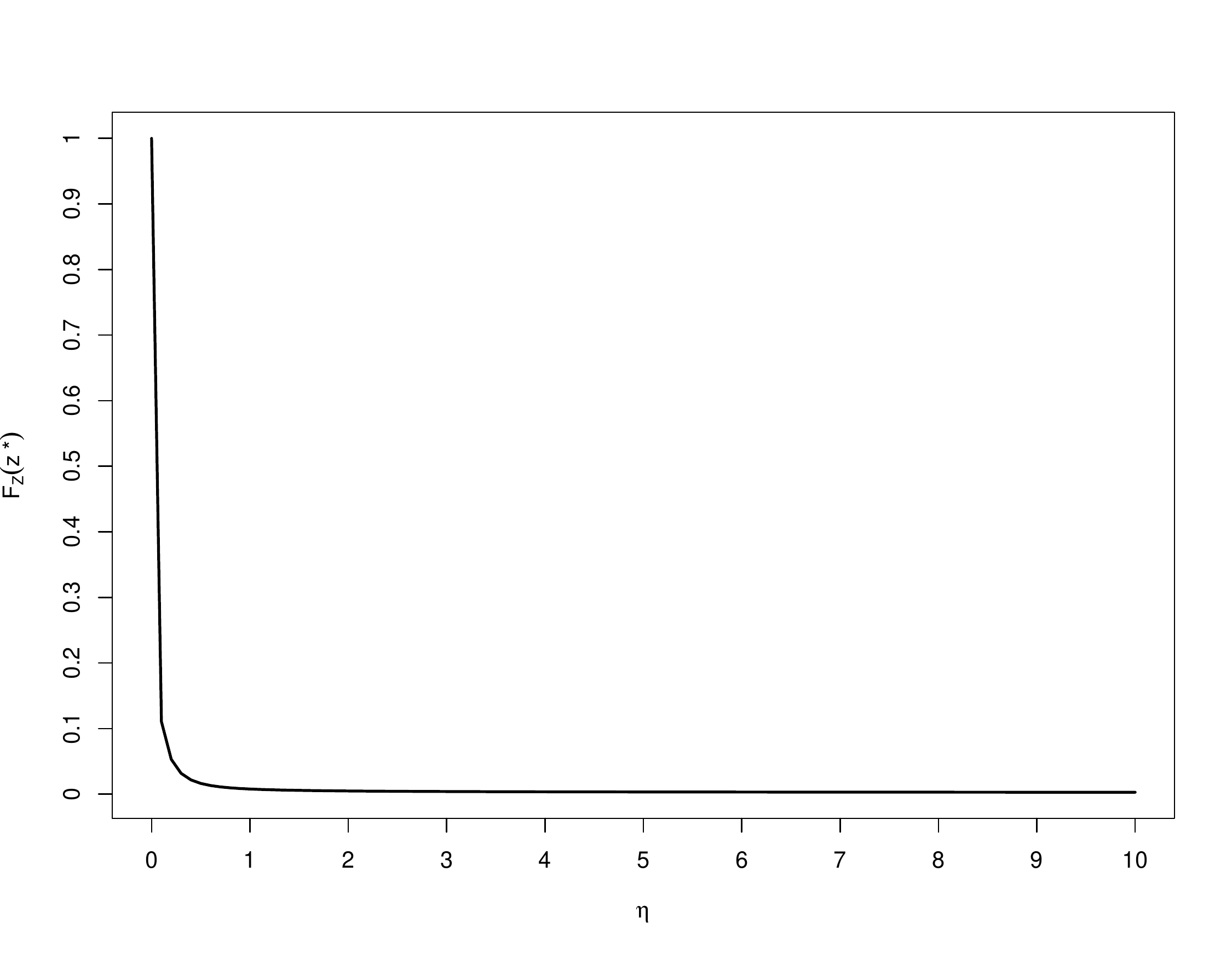}
\end{center}
\caption{The value $F_Z(z^*)$ as a function of $\eta$ for $m=20$ and $q=0.05$ under the assumption of a Clayton copula.}%
\label{Fig:Clayton-F}%
\end{figure}

\begin{figure}[ptbh]
\begin{center}%
\includegraphics[width=10.7cm]{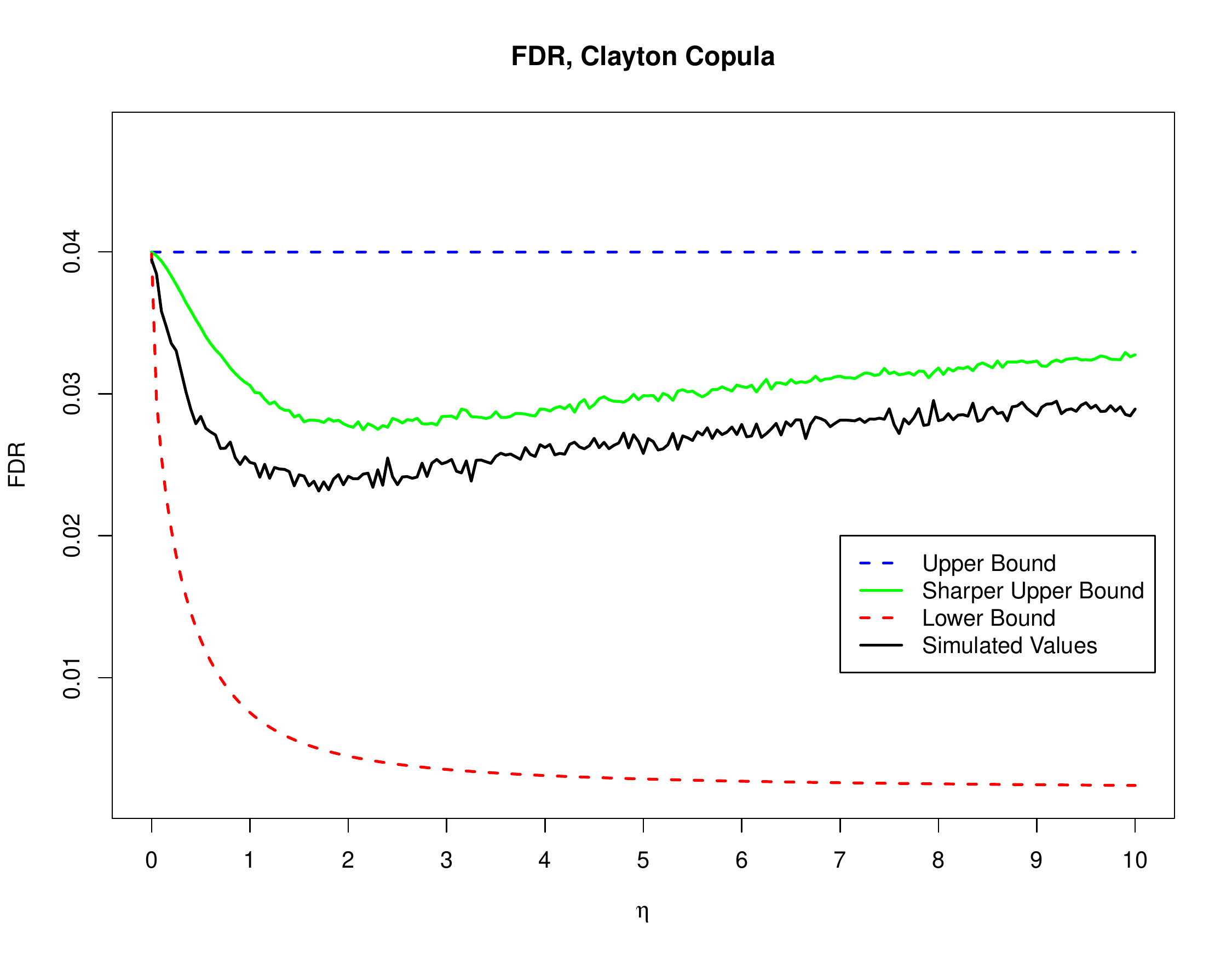}
\end{center}
\caption{Lower bound (dashed red line), upper bound (dashed blue line), the sharper upper bound (solid green line), and simulated values of the FDR of $\vp^{LSU}$ (solid black line) as functions of $\eta$ for a Clayton copula. We put $m=20$, $q=0.05$, and $m_0=16$. Simulated values are based on $10^5$ independent pseudo realizations of $Z$.}%
\label{Fig:Clayton-FDR}%
\end{figure}

The generator of the Clayton copula is given by
\begin{equation}\label{clayton_copula}
\psi(t)=(1+\eta t)^{-1 / \eta}, \quad \eta \in (0, \infty),
\end{equation}
leading to $\psi^{-1}(x)= \left(x^{-\eta}-1\right) / \eta$ and to the probability density function (pdf)
\begin{equation}\label{clayton-f_z}
f_Z(z)=\frac{1}{\eta} f_{\Gamma\left(1/\eta,1\right)}\left(z / \eta\right)=\frac{1}{\Gamma \left(1/\eta\right)}\eta^{-1/\eta} z^{1/\eta-1} \exp\left(-z / \eta\right)
\end{equation}
of $Z$, where $\Gamma$ denotes Euler's gamma function and $f_{\Gamma\left(\alpha,\beta\right)}$ the pdf of the gamma distribution with shape parameter $\alpha \in (0,\infty)$ and scale parameter $\beta \in (0,\infty)$.
For the Clayton copula, $z^*$ is given by
\begin{equation}\label{clayton-z-star}
z^*=\frac{\log m} {\eta^{-1}\left(\left(q/m\right)^{-\eta}-q^{-\eta}\right)}=\frac{\eta \log m}{\left(m / q \right)^{\eta}-\left(1 / q \right)^{\eta}}.
\end{equation}
In Figure \ref{Fig:Clayton-F}, we plot $F_Z(z^*)$ as a function of $\eta$ for $m=20$ and $q=0.05$. It is worth mentioning that the Clayton copula converges to the independence copula
for $\eta \rightarrow 0$. In this case we get $z^* \rightarrow 1$ and $f_Z(z^*)$ tends to the Dirac delta function concentrated in $1$. As a result, we observe that $F_Z(z^*)\rightarrow 1$ as $\eta \rightarrow 0$ and the FDR of $\vp^{LSU}$ approaches $m_0 q / m$. As $\eta$ increases, $F_Z(z^*)$ steeply decreases and takes values very close to zero for large values of $\eta$. Consequently, it is expected that the FDR of $\vp^{LSU}$ is close to $m_0 q / m$ for large values of $\eta$, too. For $\eta$ of moderate size, however, the FDR of $\vp^{LSU}$ can be much smaller than $m_0 q / m$. This is shown in Figure \ref{Fig:Clayton-FDR} below and discussed in detail there.

The quantity $\gamma_{min}$ for the Clayton copula is calculated by
\begin{eqnarray}
\gamma_{min}&=&1- {\eta}^{-1} \int_0^{z^*} \frac{\exp\left(-z \psi^{-1}\left(q / m \right)\right)}{q/m}  f_{\Gamma\left(1/\eta,1\right)}\left(z / \eta\right) dz + \nonumber\\
& &\;\;\;\; {\eta}^{-1} \int_0^{z^*}\frac{\exp\left(-z \psi^{-1}\left(q\right)\right)}{q} f_{\Gamma\left(1/\eta,1\right)}\left(z / \eta\right)dz \nonumber\\
&=& 1-I_1^C+I_2^C, \label{clayton_gamma_min}
\end{eqnarray}
where
\begin{eqnarray*}
I_1^C&=&\frac{\eta^{-1/\eta}}{\Gamma\left(1/\eta\right)}\frac{m}{q} \int_0^{z^*} z^{1/\eta-1} \exp\left(-\frac{z}{\eta}\left(\left(\frac{m}{q}\right)^\eta-1\right)-\frac{z}{\eta}\right) dz\\
&=&\frac{\eta^{-1/\eta}}{\Gamma\left(1/\eta\right)}\frac{m}{q} \int_0^{z^*} z^{1/\eta-1} \exp\left(-\frac{z}{\eta}\left(\frac{m}{q}\right)^\eta\right) dz\\
&=& F_{\Gamma\left(1/\eta,\eta^{-1}(m/q)^\eta\right)}(z^*)=F_{\Gamma\left(1/\eta,1\right)}(\eta^{-1}(m/q)^\eta z^*)\\
&=&F_{\Gamma\left(1/\eta,1\right)}\left(\frac{m^\eta\ln m}{m^\eta-1}\right)
\end{eqnarray*}
and, similarly,
\begin{equation*}
I_2^C=F_{\Gamma\left(1/\eta,\eta^{-1}(1/q)^\eta\right)}(z^*)=F_{\Gamma\left(1/\eta,1\right)}(\eta^{-1}(1/q)^\eta z^*)=F_{\Gamma\left(1/\eta,1\right)}\left(\frac{\ln m}{m^\eta-1}\right).
\end{equation*}
Hence, from Theorem \ref{theorem-lower-bound} we get for all $\vt \in \Theta$ that
\begin{eqnarray*} 
FDR_{\vt, \eta}(\vp^{LSU})  \geq
\frac{m_0 q}{m}\left(1+ F_{\Gamma\left(1/\eta,1\right)}\left(\frac{\ln m}{m^\eta-1}\right)-F_{\Gamma\left(1/\eta,1\right)}\left(\frac{m^\eta\ln m}{m^\eta-1}\right) \right).
\end{eqnarray*} 

Next, we consider the sharper upper bound for the FDR in the case of Clayton copulae in detail. As outlined in the discussion around Theorem \ref{coro-sharp}, a so-called Dirac-uniform configuration (cf., e.\ g., \cite{Sinica2013} and references therein) is assumed for $\bP$ in case of $m_0 < m$. Namely, the $p$-values $(P_i: i \in I_1(\vt))$ are assumed to be $\mathbb{P}_\vt$-almost surely equal to $0$. Under assumptions (I1)-(I2) from Theorem \ref{das-FDR-theorem}, Dirac-uniform configurations are least favorable (provide upper bounds) for the FDR of $\vp^{LSU}$, see \cite{benyek:2001}. In the case of dependent $p$-values, such general results are yet lacking, but it is assumed that Dirac-uniform configurations yield upper FDR bounds
for $\vp^{LSU}$ also under dependency, at least for large $m$ (cf., e.\ g., \cite{FiDi4711}). 

Under a Dirac-uniform configuration, the sharper upper bound for the FDR of $\vp^{LSU}$ is expressed by (see Theorem \ref{coro-sharp})

\begin{eqnarray}\label{sharper_upper_bound-DU}
b(m, \vt, \eta) \!\! &=& \!\! \frac{q}{m} \sum_{i=1}^{m_0} \sum_{k=m_1+1}^{m-1} \E\left[ \left(\frac{\exp \left(-Z\psi^{-1}(q_{k+1})\right)}{q_{k+1}}
- \frac{\exp \left(-Z\psi^{-1}(q_{k})\right)}{q_k}\right)\right. \nonumber\\
\!\! & & \!\! \;\;\; \times \left. \left(G_k^i(Z)-G_k^i(z^*_k)\right) \mathbf{1}_{[z^*_k,\infty)}(Z) \right],
\end{eqnarray}
where $z^*_k$ is given in (\ref{z^*_k-Sec2}) and the random set $D_{\bY;k}^{(i,z)}$ the probability of which is given by $G_k^i(z)$ can under Dirac-uniform configurations
equivalently be expressed as
\begin{equation*}
D_{\bY;k}^{(i,z)} = \Big\{\exp \left(-z\psi^{-1}\left(q_{k+1}\right)\right)\le Y_{(k)}^{(i)},..., \exp \left(-z\psi^{-1}\left(q_m\right)\right)\le Y_{(m-1)}^{(i)} \Big\}\,.
\end{equation*}
The last equality follows from the fact that $Y_{(k)}^{(i)}$, ..., $Y_{(m-1)}^{(i)}$ almost surely correspond to $p$-values associated with true null hypotheses, i.\ e.,
\[F_{P_{(k)}^{(i)}}(x)=...=F_{P_{(m-1)}^{(i)}}(x)=x\,.\]
Moreover, since each of the $Y_{(k)}^{(i)}$, ..., $Y_{(m-1)}^{(i)}$ is obtained by the same isotonic transformation from the corresponding element in the sequence $P_{(k)}^{(i)}$, ..., $P_{(m-1)}^{(i)}$, we get that $Y_{(k)}^{(i)}$, ..., $Y_{(m-1)}^{(i)}$ is an increasing sequence of independent and identically UNI$[0,1]$-distributed random variables.
Hence, the probabilities $G_k^i(z)= \mathbb{P}_{\vt, \eta}(D_{\bY;k}^{(i,z)})$ for $k \in \{m_1+1, \hdots, m-1\}$ can be calculated recursively, for instance  by making use of Bolshev's recursion (see, e.\ g., \cite{SW1986}, p. 366).

In general, Bolshev's recursion is defined in the following way. Let $0\le a_1 \le a_2 \le \hdots \le a_n \le 1$ be real constants and let $U_{(1)}\le U_{(2)}\le \hdots \le U_{(n)}$ be the order statistics of independent and identically UNI$[0,1]$-distributed random variables. We let $\bar{P}_n(a_1, \hdots, a_n)=P(a_1\le U_{(1)}, \hdots, a_n\le U_{(n)})$. Then, the probability $\bar{P}_n(a_1, \hdots, a_n)$ is calculated recursively by
\begin{equation}\label{Bolshev}
\bar{P}_n(a_1, \hdots, a_n)=1-\sum_{j=1}^n \binom{n}{j} a_j^j \bar{P}_{n-j}(a_{j+1}, \hdots, a_n).
\end{equation}
Application of (\ref{Bolshev}) with $n=m_0-1$ and
\[
a_j=\left\{
\begin{array}{ccc}
0& \text{for} & j \in \{1, \hdots, k-m_1-1\}\\
 \exp \left(-z\psi^{-1}\left(q_{j+m_1+1}\right)\right) & \text{for} & j \in \{k-m_1, \hdots, m_0-1\}\\
\end{array}
\right.
\]
for $k\in\{m_1+1, \hdots, m-1\}$ as well as numerical integration with respect to the distribution of  $Z$ over $[z_k^*,\infty]$ lead to a numerical approximation of the sharper upper bound for the FDR of $\vp^{LSU}$
under Dirac-uniform configurations.

In Figure \ref{Fig:Clayton-FDR} we present the lower bound (dashed red line), the upper bound (dashed blue line), the sharper upper bound (solid green line), and the simulated values of the FDR of $\vp^{LSU}$ (solid black line) as a function of the parameter $\eta$ of a Clayton copula. We put $m=20$, $q=0.05$, and $m_0=16$. The $p$-values which correspond to the false null hypotheses  have been set to zero. The simulated values are obtained by using $10^5$ independent repetitions. We observe that the FDR of $\vp^{LSU}$ starts at $m_0 q / m =0.04$
for $\eta = 0$ and decreases to a minimum of approximately $0.023$ at $\eta \approx 1.7$. This value is much smaller than the nominal level $q$, offering some room for improvement of $\vp^{LSU}$ for a broad range of values of $\eta$. After reaching its minimum, the FDR of $\vp^{LSU}$ increases and tends to $0.04$ as $\eta$ increases. This behavior of the FDR of $\vp^{LSU}$ is as expected from the values of $F_Z(z^*)$, as discussed around Figure \ref{Fig:Clayton-F}.

\begin{figure}[ptbh]
\begin{center}%
\includegraphics[width=10.7cm]{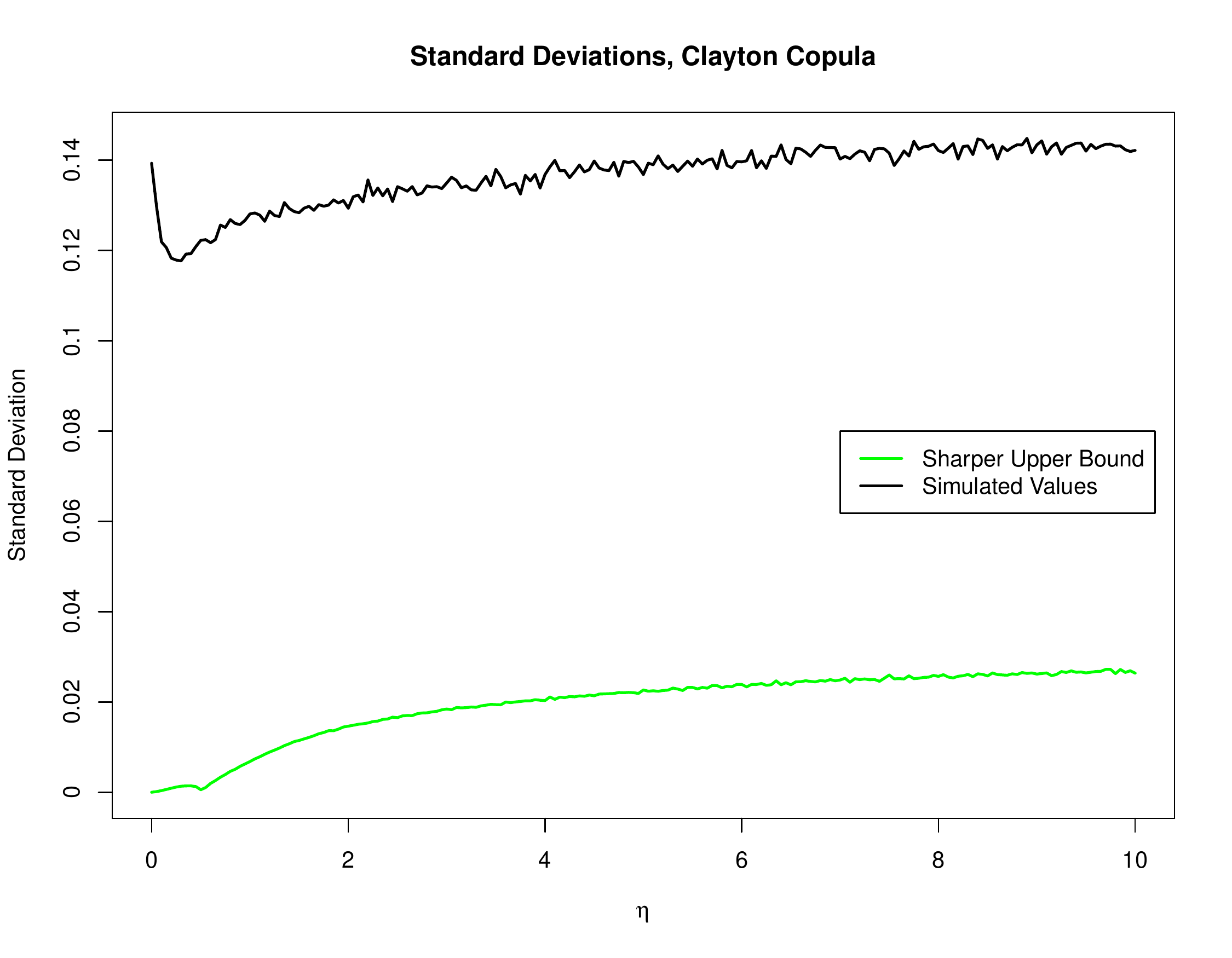}
\end{center}
\caption{Empirical standard deviations of the sharper upper bound (solid green line), and of $\FDP_{\vt, \eta}(\vp^{LSU})$ (solid black line) as functions of the parameter $\eta$ of a Clayton copula. We put $m=20$, $q=0.05$, and $m_0=16$. Simulated values are based on $10^5$ independent pseudo realizations of $Z$.}%
\label{Fig:Clayton-FDR-sd}%
\end{figure}

In contrast to the "classical" upper bound, the sharper upper bound reproduces the behavior of the simulated FDR values very well. It provides a good approximation of the true values of the FDR of $\vp^{LSU}$ for all considered values of $\eta$. In particular, it is much smaller than the "classical" upper bound for moderate values of $\eta$. Consequently, application of the sharper upper bound can be used to improve the power of the multiple testing procedure by adjusting the nominal value of $q$ depending on $\eta$. If $\eta$ is unknown, we propose techniques for pre-estimating it in Section \ref{sec5}. It is also remarkable 
that the difference between the sharper upper bound and the corresponding simulated FDR-values is not large. In contrast, the empirical standard deviations of the sharper upper bound (over repeated simulations) are about five times smaller than the 
corresponding ones for the simulated values of the FDP of $\vp^{LSU}$ (see Figure \ref{Fig:Clayton-FDR-sd}). While these standard deviations are always smaller than $0.028$ for the sharper upper bound, they are around $0.14$ for almost all of the considered values of $\eta$ in case of the simulated FDP-values. Finally, we note that the lower bound seems not to be informative in this particular model class. It is close to zero even for moderate values of $\eta$.

\subsection{Gumbel Copula}

The generator of the Gumbel copula is given by
\begin{equation}\label{gumbel_copula}
\psi(x)=\exp\left(-x^{1 / \eta}\right), \quad \eta  \geq 1,
\end{equation}
which leads to $\psi^{-1}(x)=\left(- \ln x \right)^\eta$ and a 
stochastic representation
\begin{equation} \label{Z-Gumbel}
Z\dis \left(\cos \left(\frac{\pi}{2\eta} \right)\right)^\eta Z_0, \quad \eta >1,
\end{equation}
for $Z$, where the random variable $Z_0$ has a stable distribution with index of stability $1/\eta$ and unit skewness. The cdf of $Z_0$ is given by (cf. \cite{Chambers1976}, p. 341)
\begin{eqnarray*}
F_{Z_0}(z) &= &\frac{1}{\pi}\int_0^\pi \exp \left( -z^{-1/(\eta-1)} a(v)\right) dv \quad \text{with}\\
a(v) &= &\frac{\sin\left((1-\eta)v/\eta\right) (\sin (v/\eta))^{1/(\eta-1)}}{(\sin v)^{\eta/(\eta-1)}},
\quad v \in (0,\pi).
\end{eqnarray*}
Although \eqref{Z-Gumbel} in connection with $F_{Z_0}$ characterizes the distribution of $Z$ completely, the integral representation of $F_{Z_0}$ may induce numerical issues with respect to implementation. Somewhat more convenient from this perspective is the following result. Namely, \cite{Kanter1975} obtained a stochastic representation of $Z_0$, given by
\begin{equation}\label{stable_stoch_pres}
Z_0=(a(U)/W)^{\eta-1},
\end{equation}
where $U$ and $W$ are stochastically independent, $W$ is standard exponentially distributed and $U \sim \text{UNI}(0,\pi)$. We used \eqref{stable_stoch_pres} for simulating $Z_0$ and, consequently, $Z$.

\begin{figure}[htbp]
\begin{center}%
\includegraphics[width=10.7cm]{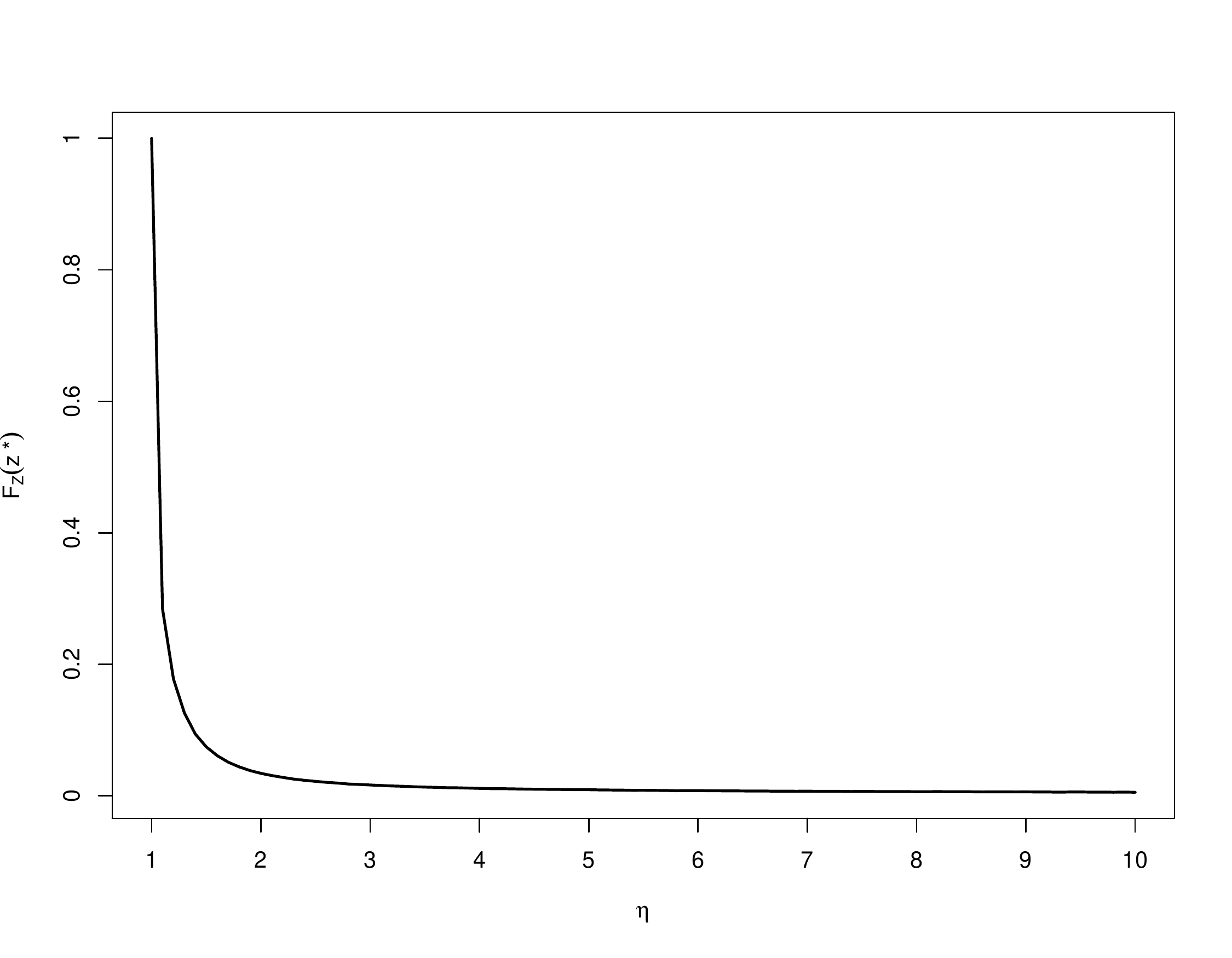}
\end{center}
\caption{The value $F_Z(z^*)$ as a function of $\eta$ for $m=20$ and $q=0.05$ under the assumption of a Gumbel copula. The graph was obtained via simulations by generating $10^6$ independent pseudo realizations of $Z$
according to \eqref{Z-Gumbel} and \eqref{stable_stoch_pres}.}%
\label{Fig:Gumbel-F}%
\end{figure}

For the Gumbel copula we get
\begin{equation}\label{gumbel-z-star}
z^*=\frac{\ln m}{\left(-\ln \frac{q}{m}\right)^\eta-\left(-\ln q \right)^\eta}=
\frac{\ln m}{\left(\ln \frac{m}{q}\right)^\eta-\left(\ln \frac{1}{q} \right)^\eta}\,.
\end{equation}
In Figure \ref{Fig:Gumbel-F}, we plot $F_Z(z^*)$ as a function of $\eta$ for $m=20$ and $q=0.05$. A similar behavior as in the case of the Clayton copula is present. If $\eta=1$ then the Gumbel copula coincides with the independence copula. Hence, $F_Z(z^*)=1$ and, consequently, the FDR of $\vp^{LSU}$ is equal to $m_0 q / m$ in this case. As $\eta$ increases, $F_Z(z^*)$ decreases and it approaches $0$ for larger values of $\eta$. Hence, $\FDR_{\vt, \eta}(\varphi^{LSU})$ tends to $m_0 q / m$ as $\eta$ becomes considerably large. For moderate values of $\eta$, $\FDR_{\vt, \eta}(\varphi^{LSU})$ can again be much smaller than $m_0 q / m$, in analogy to the situation in models with Clayton copulae.

Recall from \eqref{FDR_ineq2} that
\begin{eqnarray}
\gamma_{min}&=&1-\mathbb{E}\left[g_1(Z) \mathbf{1}_{[0,z^*]}(Z) \right], \label{gumbel_gamma_min}
\end{eqnarray}
where $g_1(Z) = g\left(\psi^{-1}(q/m)|Z\right)-g\left(\psi^{-1}(q)|Z\right)$. For the Gumbel copula, we obtain
\begin{eqnarray*}
g_1(Z)&=&\frac{\exp\left(-Z \psi^{-1}\left(q/m\right)\right)}{q/m} - \frac{\exp\left(-Z \psi^{-1}\left(q\right)\right)}{q}\\
&=&  \frac{\exp\left(-Z\left(\ln \frac{m}{q}\right)^\eta\right)}{q/m} - \frac{\exp\left(-Z \left(\ln \frac{1}{q}\right)^\eta\right)}{q}.
\end{eqnarray*}
The expectation in \eqref{gumbel_gamma_min} cannot be calculated analytically. However, it can easily be  approximated with Monte Carlo simulations by applying the stochastic representations \eqref{Z-Gumbel} and \eqref{stable_stoch_pres} for any fixed $\eta>1$. This leads to a numerical value on the left-hand side of the chain of inequalities

\begin{equation}\label{gumbel_FDR}
\frac{m_0 q}{m}\left(1-\mathbb{E}\left[g_1(Z) \mathbf{1}_{[0,z^*]}(Z) \right]\right) \le \FDR_{\vt, \eta}(\varphi^{LSU}) \le \frac{m_0 q}{m}.
\end{equation}

\begin{figure}[htbp]
\begin{center}%
\includegraphics[width=10.7cm]{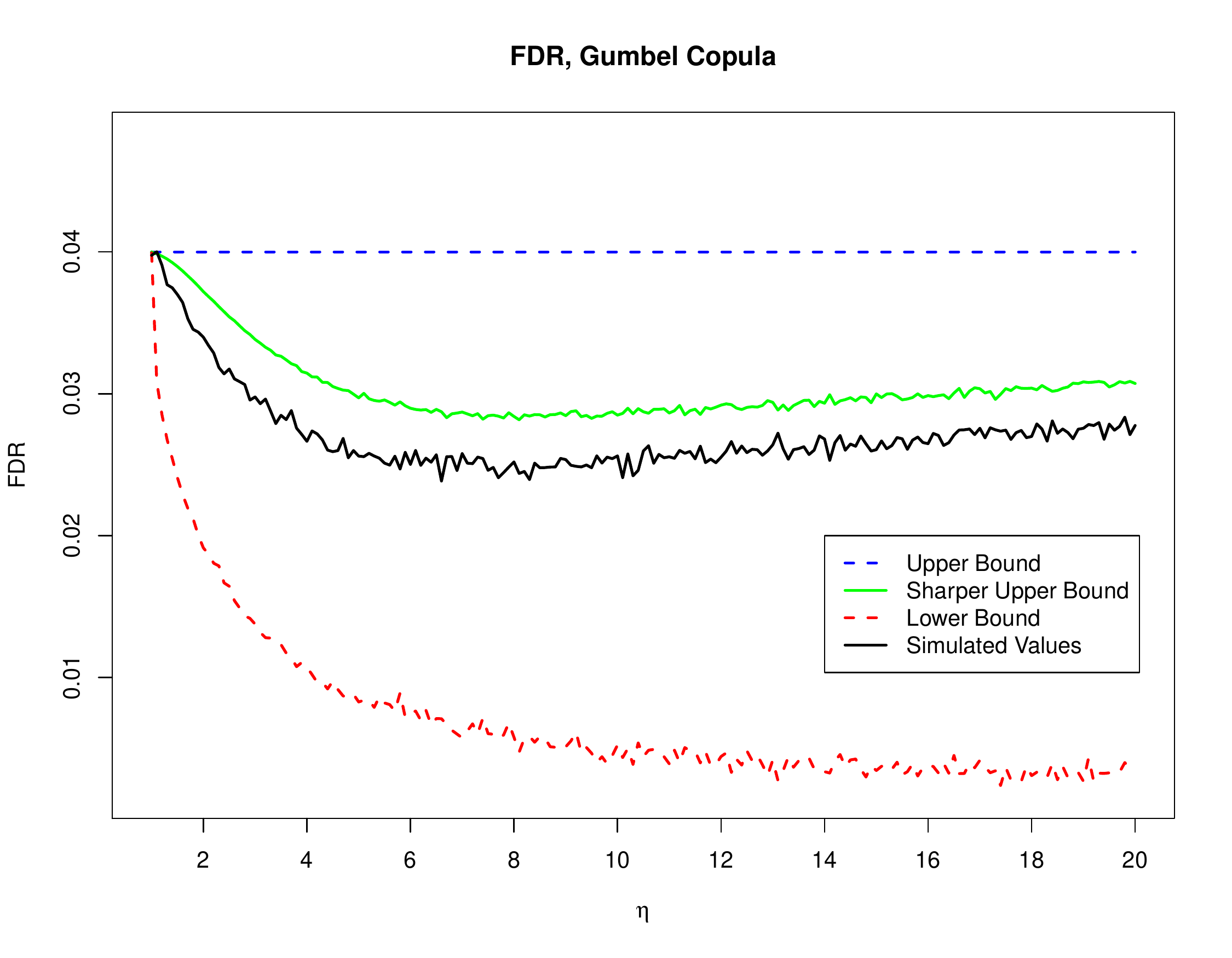}
\end{center}
\caption{Lower bound (dashed red line), upper bound (dashed blue line), the sharper upper bound (solid green line), and simulated values of the FDR of $\vp^{LSU}$ (solid black line) as functions of the parameter $\eta$ of a Gumbel copula. We put $m=20$, $q=0.05$, and $m_0=16$. Simulated values are based on $10^5$ independent pseudo realizations of $Z$.}%
\label{Fig:Gumbel-FDR}%
\end{figure}

The sharper upper bound from Theorem \ref{coro-sharp} can be calculated by using Bolshev's recursion similarly to the discussion around \eqref{Bolshev}, but here with $\psi$ as in \eqref{gumbel_copula}. Figure \ref{Fig:Gumbel-FDR} displays the lower bound (dashed red line), the upper bound (dashed blue line), the sharper upper bound (solid green line), and simulated values of $\FDR_{\vt, \eta}(\varphi^{LSU})$ (solid black line) as functions of $\eta$.
Again, we chose $m=20$, $q=0.05$, and $m_0= 16$. The $p$-values corresponding to the false null hypotheses were all set to zero, as in the case of Clayton copulae. The simulated values were obtained by generating $10^5$ independent pseudo realizations of $Z$.

\begin{figure}[htbp]
\begin{center}%
\includegraphics[width=10.7cm]{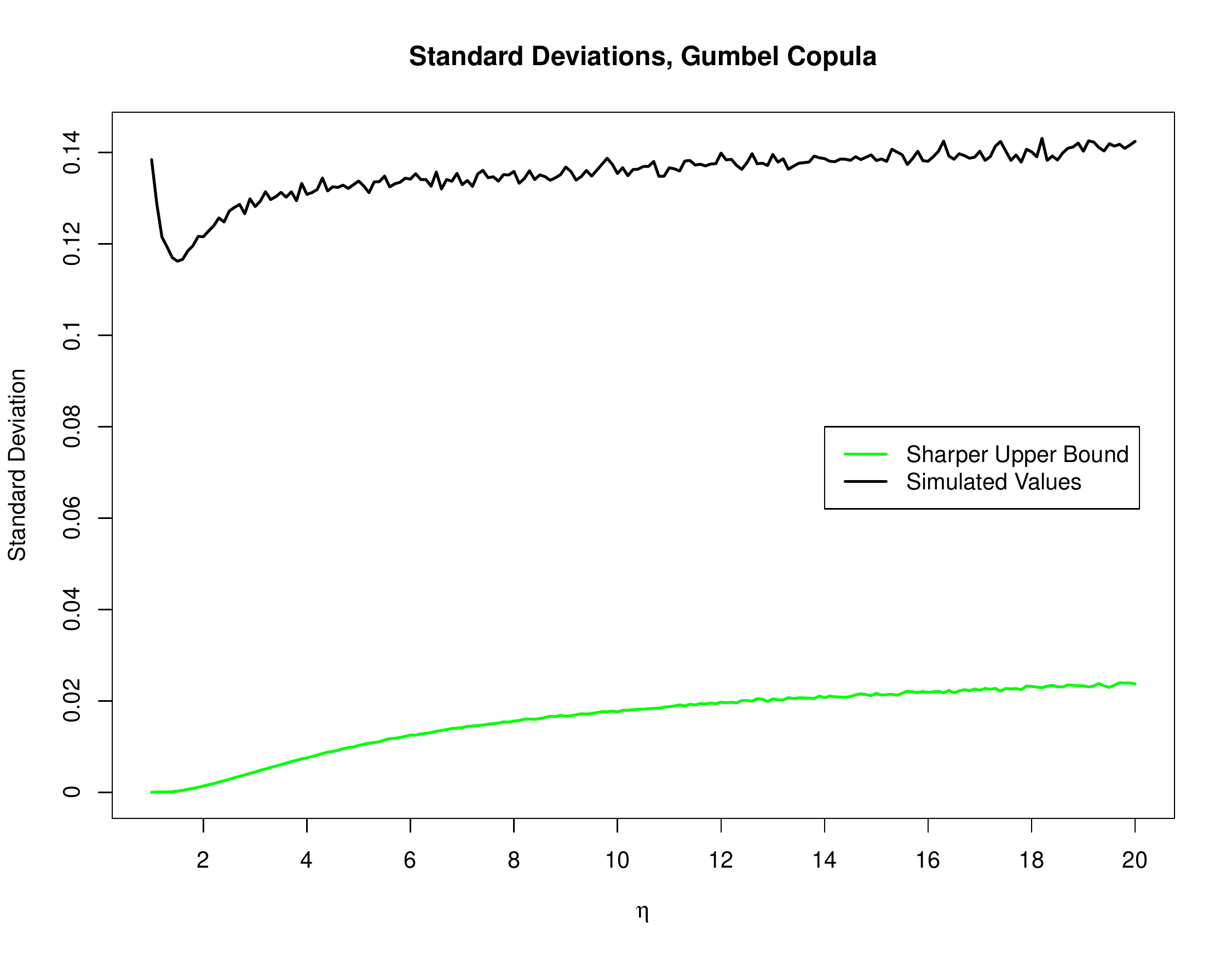}
\end{center}
\caption{Empirical standard deviations of the sharper upper bound (solid green line), and of $\FDP_{\vt, \eta}(\vp^{LSU})$ (solid black line) as functions of the parameter $\eta$ of a Gumbel copula. We put $m=20$, $q=0.05$, and $m_0=16$. Simulated values are based on $10^5$ independent pseudo realizations of $Z$.}%
\label{Fig:Gumbel-FDR-sd}%
\end{figure}

Similarly to the case of the Clayton copula, the curve of simulated FDR values has a $U$-shape. It starts at $m_0 q / m =0.04$ and drops to its minimum of approximately $0.024$ for values of $\eta$ around $6.6$.
For such values of $\eta$, the black curve is considerably below the classical upper bound of $0.04$. In contrast, the sharper upper bound gives a much tighter approximation of the simulated FDR values in such cases
and reproduces the $U$-shape over the entire range of values for the parameter $\eta$ of the Gumbel copula.
As a result, its application can be used to improve power by adjusting the nominal value of $q$ and thereby increasing the probability to detect false null hypotheses. Moreover, as in the case of Clayton copulae, the empirical standard deviations of the sharper upper bound are much smaller than those of the simulated values of the FDP (see Figure \ref{Fig:Gumbel-FDR-sd}). 
The lower bound from \eqref{gumbel_FDR} (corresponding to the dashed red curve in Figure \ref{Fig:Gumbel-FDR}) has been obtained by approximating the expectation in \eqref{gumbel_gamma_min} via simulations. As in the case of the Clayton copula, the lower bound is not too informative for the model class that we have considered here (Dirac-uniform configurations).

\section{Empirical copula calibration} \label{sec5}
In the previous section we studied the influence of the copula parameter $\eta$ on the FDR of $\vp^{LSU}$ under several parametric families of Archimedean copulae. It turned out that adapting $\vp^{LSU}$ to the degree of dependency in the data by adjusting the nominal value of $q$ based on the sharper upper bound from Theorem
\ref{coro-sharp} is a promising idea, because the unadjusted procedure may lead to a considerable non-exhaustion of $q$, cf.\ Figures \ref{Fig:Clayton-FDR} and \ref{Fig:Gumbel-FDR}. Due to the decision rule of a step-up test, this also entails suboptimal power properties of $\vp^{LSU}$ when applied "as is" to models with Archimedean $p$-value copulae.

In practice, however, often the copula parameter itself is an unknown quantity. Hence, the outlined adaptation of $q$ typically requires some kind of pre-estimation of $\eta$ before multiple testing is performed.
Although this is not in the main focus of the present work, we therefore outline possibilities for estimating $\eta$ and for quantifying the uncertainty of the estimation in this section.

One class of procedures relies on resampling, namely via the parametric bootstrap or via permutation techniques if $H_1, \hdots, H_m$ correspond to marginal two-sample problems. \cite{Pollard2004} provided an extensive comparison of both approaches and argued that the permutation method 
reproduces the correct null distribution only under some conditions. However, if these conditions are met, the permutation approach is often superior to bootstrapping (see also \cite{westfallyoung} and \cite{meinshausen-optimal}). Furthermore, it is essential to keep in mind that both bootstrap and permutation-based methods estimate the distribution of the vector $\bP$ under the global null hypothesis $H_0$. Hence, the assumption that $\eta$ does not depend on $\vt$ is an essential prerequisite for the applicability of such resampling methods for estimating $\eta$. Notice that the latter assumption is an informal description of the "subset pivotality" condition introduced by \cite{westfallyoung}. The resampling methods developed by \cite{dudoit2008} can dispense with subset pivotality in special model classes, but for the particular task of estimating the copula parameter this assumption seems indispensable.

Estimation of $\eta$ and uncertainty
quantification of the estimation based on resampling is generally performed by applying a suitable
estimator $\hat{\eta}$ to the re- (pseudo) samples. In the context of Archimedean copulae the two most widely applied estimation procedures are the maximum likelihood method (see, e.\ g. \cite{Joe2005}, \cite{Hofert2012}) and the method of moments (referred to as "realized copula" approach by \cite{fengler2012}).

\cite{Hofert2012} considered the estimation of the parameter of an Archimedean copula with known margins by the maximum likelihood approach. To this end, they derived analytic expressions for the derivatives of the copula generator for several families of Archimedean copulae, as well as formulas for the corresponding score functions. Using these results and assuming a regular model, an elliptical asymptotic confidence region for the copula parameter $\eta$ can be obtained by applying general limit theorems for maximum likelihood estimators (see \cite{Hofert2012} for details and the calculations for different types of Archimedean copulae).

In the context of the method of moments, Kendall's tau is often considered. For a bivariate Archimedean copula
with generator $\psi$ of marginally UNI$[0, 1]$-distributed variates $P_1$ and $P_2$, it is given by
\begin{eqnarray}
\tau_{P_1,P_2}&=&4\int_0^1\int_0^1 F_{(P_1,P_2)}(u,v)dF_{(P_1,P_2)}(u,v)-1 \nonumber \\
&=& 1-4\int_0^{\psi^{-1}(0)} t [\psi^\prime(t)]^2 dt, \label{tau-psi}
\end{eqnarray}
cf. \cite{McNeil2009}.

The right-hand side of \eqref{tau-psi} can analytically be calculated for some families of Archimedean copulae. For instance, for a Clayton copula with parameter $\eta$ it is given by $\tau(\eta)=\eta / (2+\eta)$, while it is equal to $\tau(\eta)=(\eta-1) / \eta$ for a Gumbel copula with parameter $\eta$ (see \cite{nelsen2006}, p. 163-164). Based on such moment equations, \cite{fengler2012} suggested the "realized copula" method for empirical calibration
of a one-dimensional parameter $\eta$ of an $m$-variate Archimedean copula. The method considers all  $m(m -1)/2$ distinct pairs of the
$m$ underlying random variables, replaces the population versions of $\tau(\eta)$ by the corresponding sample analogues, and finally aggregates the resulting $m(m-1)/2$
estimates in an appropriate manner. More specifically, consider the functions $g_{ij}(\eta) = \hat{\tau}_{ij}-\tau(\eta)$ for $1 \le i < j \le m$ and define $\mathbf{q}(\eta) = (g_{ij}(\eta): 1\le i<j\le m)^\top$, where $\hat{\tau}_{ij}$ is the sample estimator of Kendall's tau (see, e.\ g., \cite{nelsen2006}, Section 5.1.1). The resulting estimator for $\eta$ is then obtained by
\begin{equation}\label{RC}
\hat{\eta} = \arg \min_{\eta} \left\{ \mathbf{q}(\eta)^\top \mathbf{W} \mathbf{q}(\eta) \right\}
\end{equation}
for an appropriate weight matrix $\mathbf{W} \in \mathbb{R}^{\binom{m}{2} \times \binom{m}{2}}$.
An application of the realized copula method to resampled $p$-values generated by permutations in the context of multiple testing for differential gene expression has been demonstrated  by \cite{stp-copulae}.
Multivariate extensions of Kendall's tau and central limit theorems for the sample versions have been derived by \cite{genest2011}. These results can be used for uncertainty quantification of the moment estimation
of $\eta$ by constructing asymptotic confidence regions.

\section{Discussion} \label{sec6}
We have derive a sharper upper bound for the FDR of $\vp^{LSU}$ in models with Archimedean copulae.  This bound can be used to prove that $\vp^{LSU}$ controls the FDR for this type of multivariate $p$-value distributions, a result which is in line with the findings of \cite{benyek:2001} and \cite{Sarkar2002}. Since certain models with $H_0$-exchangeable $p$-values fall into this class at least asymptotically (see Theorem \ref{thm-ex}), our findings complement those of \cite{FiDi4711} who investigated infinite sequences of $H_0$-exchangeable $p$-values in Gaussian models. While our general results in Section \ref{sec3} qualitatively extend the theory, our results in Section \ref{sec4} regarding Clayton and Gumbel copulae are quantitatively very much in line with the findings for Gaussian and $t$-copulae reported by \cite{FiDi4711}. Namely, over a broad class of models with dependent $p$-values, the FDR of $\vp^{LSU}$ as a function of the dependency parameter has a $U$-shape and becomes smallest for medium strength of dependency among the $p$-values. This behavior can be exploited by adjusting $q$ in order to adapt to $\eta$. We have presented an explicit adaptation scheme based on the upper bound from Theorem \ref{coro-sharp}. To the best of our knowledge, this kind of adaptation is novel to FDR theory.

It is beyond the scope of the present work to investigate which parametric class of copulae is appropriate for which kind of real-life application. Relatedly, the problem of model misspecification (i.\ e., quantification of the approximation error if the true model does not belong to the class with Archimedean $p$-value copulae and is approximated by the (in some suitable norm) closest member of this class) could not be addressed here, but is a challenging topic for future research. One particularly interesting issue in this direction is FDR control for finite sequences of $H_0$-exchangeable $p$-values.

Finally, we would like to mention that the empirical variance of the false discovery proportion was large in all our simulations, implying that the random variable $\FDP_{\vt, \eta}(\vp^{LSU})$ was not well concentrated around its expected value $\FDR_{\vt, \eta}(\vp^{LSU})$. This is a known effect for models with dependent $p$-values (see, e.\ g., \cite{FiDi4711}, \cite{Roquain-equi}, \cite{Sinica2013}) and provokes the question if FDR control is a suitable criterion under dependency at all. Maybe more stringent in dependent models is control of the false discovery exceedance rate, meaning to design a multiple test $\vp$ ensuring that
$\FDX_{\vt, \eta}(\vp) = \mathbb{P}_{\vt, \eta}(\FDP_{\vt, \eta}(\vp) > c) \leq \gamma$, for user-defined parameters $c$ and $\gamma$. In any case, practitioners should be (made) aware of the fact that controlling the FDR with $\vp^{LSU}$ does not necessarily imply that the FDP for their particular experiment is small, at least if dependencies among $P_1, \hdots, P_m$ have to be assumed as it is typically the case in applications. In contrast, the empirical standard deviations of our proposed sharper upper bound are about five times smaller than the empirical standard deviations of the simulated values of the FDP of $\vp^{LSU}$. This provides an additional (robustness) argument for the application of the results presented in Theorem \ref{coro-sharp} in practice.

\section{Proofs} \label{sec_proofs}

\subsection*{Proof of Theorem \ref{thm1}}
Following \cite{benyek:2001}, an analytic expression for the FDR of $\vp^{LSU}$ is given by
\begin{equation}\label{FDR_expr1}
\FDR_{\vt, \eta}(\vp^{LSU}) =\sum_{i=1}^{m_0} \sum_{k=1}^{m} \frac{1}{k} \mathbb{P}_{\vt, \eta}\left\{A_k^{(i)}\right\},
\end{equation}
where $A_k^{(i)}=\{P_i\le q_k \cap C_k^{(i)}\}$ denotes the event that $k$ hypotheses are rejected one of which is $H_{i}$ (a true null hypothesis) and $C_k^{(i)}$ is the event that $k-1$ hypotheses additionally to $H_{i}$ are rejected. It holds that $(C_k^{(i)}: 1 \leq k \leq m)$ are disjoint and that $\bigcup_{k=1}^m C_k^{(i)}=[0,1]^{m-1}$.

Let $D_k^{(i)}=\bigcup_{j=1}^k C_j^{(i)}$ for $k =1,\hdots,m$ denote the event that the number of rejected null hypotheses is at most $k$. In terms of $\bP^{(i)}$ introduced in Theorem \ref{coro-sharp}, the random set $D_k^{(i)}$ is given by
\begin{equation} \label{D_k-i}
D_k^{(i)} = \{q_{k+1}\le P_{(k)}^{(i)},\hdots, q_{m}\le P_{(m-1)}^{(i)} \}.
\end{equation}

Next, we prove that
\begin{eqnarray}\label{Th2_1}
\frac{\mathbb{P}_{\vt, \eta}\left(P_i\le q_k \cap D_k^{(i)}\right)}{q_k} &\le& \frac{\mathbb{P}_{\vt, \eta}\left(P_i\le q_{k+1} \cap D_k^{(i)}\right)}{q_{k+1}}\\
&-&\int_{z^*_k}^\infty \left(\frac{\exp \left(-z\psi^{-1}(q_{k+1})\right)}{q_{k+1}} - \frac{\exp \left(-z\psi^{-1}(q_{k})\right)}{q_k}\right) \times \nonumber\\
& & \hspace*{40pt} (G_k^i(z)-G_k^i(z^*_k)) dF_Z(z)\,.\nonumber
\end{eqnarray}

To this end, we consider the function $\mathbf{T}$ introduced in Theorem \ref{coro-sharp}, which transforms a possible realization of the original $p$-values $\bP$ into a realization of $\bY$ for $Z=z$, where
$\bY = (Y_1, \hdots, Y_m)^\top$ and $Z$ are as in \eqref{P-stoch_pres}.
Because each component of this multivariate transformation is a monotonically increasing function which fully covers the interval $[0,1]$, the resulting transformation bijectively transforms the set $[0,1]^m$ into itself.
Let $C_{\bY;k}^{(i,z)}$ and $D_{\bY;k}^{(i,z)}$ denote the images of the sets $C_{k}^{(i)}$ and $D_{k}^{(i)}$ under $\mathbf{T}$ for given $Z=z$. Then
\begin{itemize}
\item[(a)] $C_{\bY;k}^{(i,z)}$ are disjoint, i.\ e., $C_{\bY;k_1}^{(i,z)} \cap C_{\bY;k_2}^{(i,z)}=\emptyset$ for $1 \leq k_1 \neq k_2 \leq m$,
\item[(b)] $D_{\bY;k}^{(i,z)}=\bigcup_{j=1}^k C_{\bY;j}^{(i,z)}$,
\item[(c)] $D_{\bY;m}^{(i,z)}=\bigcup_{j=1}^m C_{\bY;j}^{(i,z)}=[0,1]^{m-1}$.
\end{itemize}
Statements (a) - (c) follow directly from the facts that each $T_j$ is a monotonically increasing function and $\textbf{T}$ is a one-to-one transformation with image equal to $[0,1]^m$. Moreover, we obtain
\begin{equation}\label{Th2_2}
D_{\bY;k}^{(i,z)} = \Big\{ \forall k \leq j \leq m-1: Y_{(j)}^{(i)} \geq \exp \left(-z\psi^{-1}\left(F_{P_{(j)}^{(i)}}(q_{j+1})\right)\right) \Big\},
\end{equation}
where $\bY^{(i)}$ is the $(m-1)$-dimensional vector obtained from $\bY=(Y_1,\hdots,Y_m)^T$ by deleting $Y_i$.
The last equality shows that $D_{\bY;k}^{(i,z_1)} \subseteq D_{\bY;k}^{(i,z_2)}$ for $z_1 \le z_2$ and, hence,
that $G_k^i$, given by $G_k^i(z) = \mathbb{P}_{\vt, \eta} \left(D_{\bY;k}^{(i,z)}\right)$,
is an increasing function in $z$.

Returning to \eqref{Th2_1}, we obtain
\begin{eqnarray}
&& \frac{\mathbb{P}_{\vt, \eta}\left(P_i\le q_{k+1} \cap D_k^{(i)}\right)}{q_{k+1}} - \frac{\mathbb{P}_{\vt, \eta}\left(P_i\le q_k \cap D_k^{(i)}\right)}{q_k} \nonumber \\
&=& \int \left(\frac{\mathbb{P}_{\vt, \eta}\left(P_i\le q_{k+1} \cap D_k^{(i)}|Z=z\right)}{q_{k+1}}
- \right. \nonumber\\
& & \hspace*{47pt} \left. \frac{\mathbb{P}_{\vt, \eta}\left(P_i\le q_k \cap D_k^{(i)}|Z=z\right)}{q_k}\right) dF_Z(z) \nonumber\\
&=& \int \left(\frac{\mathbb{P}_{\vt, \eta}\left(P_i\le q_{k+1}|Z=z\right)\mathbb{P}_{\vt, \eta}\left(D_k^{(i)}|Z=z\right)}{q_{k+1}} - \right. \nonumber\\
& & \hspace*{47pt} \left.\frac{\mathbb{P}_{\vt, \eta}\left(P_i\le q_{k}|Z=z\right)\mathbb{P}_{\vt, \eta}\left(D_k^{(i)}|Z=z\right)}{q_k}\right) dF_Z(z)\nonumber\\
&=& \int \left(\frac{\mathbb{P}_{\vt, \eta}\left(Y_i\le \exp \left(-z\psi^{-1}(q_{k+1})\right)\right)}{q_{k+1}} - \right. \nonumber\\
& & \hspace*{47pt} \left. \frac{\mathbb{P}_{\vt, \eta}\left(Y_i\le \exp \left(-z\psi^{-1}(q_{k})\right)\right)}{q_k}\right)\mathbb{P}_{\vt, \eta}\left(D_{\bY;k}^{(i,z)}\right) dF_Z(z)\nonumber\\
&=& \int \left(\frac{\exp \left(-z\psi^{-1}(q_{k+1})\right)}{q_{k+1}} - \frac{\exp \left(-z\psi^{-1}(q_{k})\right)}{q_k}\right)G_k^i(z) dF_Z(z).\label{Th2_4}
\end{eqnarray}

Next, we analyze the difference under the last integral. It holds that
\begin{eqnarray*}
&&\frac{\exp \left(-z\psi^{-1}(q_{k+1})\right)}{q_{k+1}} - \frac{\exp \left(-z\psi^{-1}(q_{k+1})\right)}{q_k}\\
&=& \exp \left(- \log q_{k+1} -z\psi^{-1}(q_{k+1})\right) - \exp \left(- \log q_{k} -z\psi^{-1}(q_{k})\right)\\
&=& \exp \left(- \log q_{k} -z\psi^{-1}(q_{k})\right) \times \\
& &\;\;\; \left( \exp \left(- \log q_{k+1}+\log q_{k} -z\psi^{-1}(q_{k+1})+z\psi^{-1}(q_{k})\right) -1 \right).
\end{eqnarray*}
The last expression is nonnegative if and only if
\[- \log q_{k+1}+\log q_{k} -z\psi^{-1}(q_{k+1})+z\psi^{-1}(q_{k}) \ge 0 \,.\]
Hence, for $z \ge z^*_k$ with $z^*_k$ given in \eqref{z^*_k-Sec2},
the function under the integral in (\ref{Th2_4}) is positive and for $z \le z^*_k$ it is negative.
Application of this result leads to
\begin{eqnarray*}
&& \frac{\mathbb{P}_{\vt, \eta}\left(P_i\le q_{k+1} \cap D_k^{(i)}\right)}{q_{k+1}} - \frac{\mathbb{P}_{\vt, \eta}\left(P_i\le q_k \cap D_k^{(i)}\right)}{q_k} \\
&=& \int_{0}^{z^*_k} \left(\frac{\exp \left(-z\psi^{-1}(q_{k+1})\right)}{q_{k+1}} - \frac{\exp \left(-z\psi^{-1}(q_{k})\right)}{q_k}\right)G_k^i(z) dF_Z(z)\\
&+& \int_{z^*_k}^\infty \left(\frac{\exp \left(-z\psi^{-1}(q_{k+1})\right)}{q_{k+1}} - \frac{\exp \left(-z\psi^{-1}(q_{k})\right)}{q_k}\right)G_k^i(z) dF_Z(z)\\
&\ge& G_k^i(z^*_k) \int_{0}^{z^*_k} \left(\frac{\exp \left(-z\psi^{-1}(q_{k+1})\right)}{q_{k+1}} - \frac{\exp \left(-z\psi^{-1}(q_{k})\right)}{q_k}\right) dF_Z(z)\\
&+&\int_{z^*_k}^\infty \left(\frac{\exp \left(-z\psi^{-1}(q_{k+1})\right)}{q_{k+1}} - \frac{\exp \left(-z\psi^{-1}(q_{k})\right)}{q_k}\right)G_k^i(z) dF_Z(z)\,.
\end{eqnarray*}

Because of
\begin{eqnarray*}
&& \int \left(\frac{\exp \left(-z\psi^{-1}(q_{k+1})\right)}{q_{k+1}} - \frac{\exp \left(-z\psi^{-1}(q_{k})\right)}{q_k}\right) dF_Z(z)\\
&=& \int \frac{\exp \left(-z\psi^{-1}(q_{k+1})\right)}{q_{k+1}} dF_Z(z) - \int \frac{\exp \left(-z\psi^{-1}(q_{k})\right)}{q_k} dF_Z(z)\\
&=& \frac{\psi \left(\psi^{-1}(q_{k+1})\right)}{q_{k+1}}- \frac{\psi\left(\psi^{-1}(q_{k})\right)}{q_k}=\frac{q_{k+1}}{q_{k+1}}-\frac{q_k}{q_k}=0
\end{eqnarray*}
we get
\begin{eqnarray*}
&&\int_{0}^{z^*_k} \left(\frac{\exp \left(-z\psi^{-1}(q_{k+1})\right)}{q_{k+1}} - \frac{\exp \left(-z\psi^{-1}(q_{k})\right)}{q_k}\right) dF_Z(z)\\
&=&-\int_{z^*_k}^\infty \left(\frac{\exp \left(-z\psi^{-1}(q_{k+1})\right)}{q_{k+1}} - \frac{\exp \left(-z\psi^{-1}(q_{k})\right)}{q_k}\right) dF_Z(z)
\end{eqnarray*}
and, consequently,
\begin{eqnarray}
&& \frac{\mathbb{P}_{\vt, \eta}\left(P_i\le q_{k+1} \cap D_k^{(i)}\right)}{q_{k+1}} - \frac{\mathbb{P}_{\vt, \eta}\left(P_i\le q_k \cap D_k^{(i)}\right)}{q_k} \nonumber\\
&\ge&- G_k^i(z^*_k)\int_{z^*_k}^\infty \left(\frac{\exp \left(-z\psi^{-1}(q_{k+1})\right)}{q_{k+1}} - \frac{\exp \left(-z\psi^{-1}(q_{k})\right)}{q_k}\right) dF_Z(z)\nonumber\\
&+&\int_{z^*_k}^\infty \left(\frac{\exp \left(-z\psi^{-1}(q_{k+1})\right)}{q_{k+1}} - \frac{\exp \left(-z\psi^{-1}(q_{k})\right)}{q_k}\right)G_k^i(z) dF_Z(z)\nonumber\\
&=&\int_{z^*_k}^\infty \left(\frac{\exp \left(-z\psi^{-1}(q_{k+1})\right)}{q_{k+1}} - \frac{\exp \left(-z\psi^{-1}(q_{k})\right)}{q_k}\right) \times \label{lower_bound}\\
& & \hspace*{40pt} (G_k^i(z)-G_k^i(z^*_k)) dF_Z(z) ,\nonumber
\end{eqnarray}
which is obviously positive since both the differences under the integral in \eqref{lower_bound} are positive. This completes the proof of \eqref{Th2_1}.

Using (\ref{Th2_1}), we get for all $1 \le k \le m-1$ that
\begin{eqnarray*}
&&\frac{\mathbb{P}_{\vt, \eta}\left(P_i\le q_k \cap D_k^{(i)}\right)}{q_k}+\frac{\mathbb{P}_{\vt, \eta}\left(P_i\le q_{k+1} \cap C_{k+1}^{(i)}\right)}{q_{k+1}}\\
&\le&\frac{\mathbb{P}_{\vt, \eta}\left(P_i\le q_{k+1} \cap D_k^{(i)}\right)}{q_{k+1}}+\frac{\mathbb{P}_{\vt, \eta}\left(P_i\le q_{k+1} \cap C_{k+1}^{(i)}\right)}{q_{k+1}}\\
&-&\int_{z^*_k}^\infty \left(\frac{\exp \left(-z\psi^{-1}(q_{k+1})\right)}{q_{k+1}} - \frac{\exp \left(-z\psi^{-1}(q_{k})\right)}{q_k}\right) \times\\
& & \hspace*{40pt} (G_k^i(z)-G_k^i(z^*_k)) dF_Z(z)\\
&=&\frac{\mathbb{P}_{\vt, \eta}\left(P_i\le q_{k+1} \cap D_{k+1}^{(i)}\right)}{q_{k+1}}\\
&-&\int_{z^*_k}^\infty \left(\frac{\exp \left(-z\psi^{-1}(q_{k+1})\right)}{q_{k+1}} - \frac{\exp \left(-z\psi^{-1}(q_{k})\right)}{q_k}\right) \times \\
& & \hspace*{40pt} (G_k^i(z)-G_k^i(z^*_k)) dF_Z(z)
\end{eqnarray*}
and, consequently, starting with $D_1^{(i)}=C_1^{(i)}$ and proceeding step-by-step for all $k \le m-1$, we obtain
\begin{eqnarray*}
&&\sum_{k=1}^{m} \frac{\mathbb{P}_{\vt, \eta}\left\{P_i\le q_{k+1} \cap C_k^{(i)}\right\}}{q_k}\le \frac{\mathbb{P}_{\vt, \eta}\left\{P_i\le q_{m} \cap D_m^{(i)}\right\}}{q_m}\\
&-&\sum_{k=1}^{m-1}\int_{z^*_k}^\infty \left(\frac{\exp \left(-z\psi^{-1}(q_{k+1})\right)}{q_{k+1}} - \frac{\exp \left(-z\psi^{-1}(q_{k})\right)}{q_k}\right) \times \\
& & \hspace*{40pt} (G_k^i(z)-G_k^i(z^*_k)) dF_Z(z)\\
&=&1-\sum_{k=1}^{m-1}\int_{z^*_k}^\infty \left(\frac{\exp \left(-z\psi^{-1}(q_{k+1})\right)}{q_{k+1}} - \frac{\exp \left(-z\psi^{-1}(q_{k})\right)}{q_k}\right) \times \\
& & \hspace*{40pt} (G_k^i(z)-G_k^i(z^*_k)) dF_Z(z).
\end{eqnarray*}
Hence,
\begin{eqnarray*}
&&\FDR_{\vt, \eta}(\vp^{LSU})= \sum_{i=1}^{m_0} \sum_{k=1}^{m} \frac{1}{k} \mathbb{P}_{\vt, \eta}\left\{A_k^{(i)}\right\}\\
&=& \sum_{i=1}^{m_0} \frac{q}{m} \sum_{k=1}^{m} \frac{\mathbb{P}_{\vt, \eta}\left\{P_i\le q_{k+1} \cap C_k^{(i)}\right\}}{q_k}\\
&\le&\sum_{i=1}^{m_0} \frac{q}{m}-\sum_{i=1}^{m_0} \frac{q}{m} \sum_{k=1}^{m-1}\int_{z^*_k}^\infty \left(\frac{\exp \left(-z\psi^{-1}(q_{k+1})\right)}{q_{k+1}} - \frac{\exp \left(-z\psi^{-1}(q_{k})\right)}{q_k}\right) \times \\
& & \hspace*{40pt} (G_k^i(z)-G_k^i(z^*_k)) dF_Z(z)\\
&=&\frac{m_0}{m}q-b(m, \vt, \psi)\,,
\end{eqnarray*}
where $b(m, \vt, \psi)$ is defined in Theorem \ref{coro-sharp}. This completes the proof of the theorem.\qed

\subsection*{Proof of Theorem \ref{theorem-lower-bound}}
Straightforward calculation yields
\begin{eqnarray*}
\FDR_{\vt, \eta}(\vp^{LSU})&=& \sum_{i=1}^{m_0} \sum_{k=1}^{m} \frac{1}{k}
\int \mathbb{P}_{\vt, \eta}\left\{ A_k^{(i)} |Z=z\right\}d F_Z(z)\\
&=&\sum_{i=1}^{m_0} \sum_{k=1}^{m} \frac{1}{k}
\int \mathbb{P}_{\vt, \eta}\left\{P_i\le q_k|Z=z\right\} \times \\
& & \hspace*{68pt} \mathbb{P}_{\vt, \eta}\left\{ C_k^{(i)}|Z=z\right\}d F_Z(z)\\
&=&\sum_{i=1}^{m_0} \frac{q}{m} \sum_{k=1}^{m}
\int \frac{\mathbb{P}_{\vt, \eta}\left\{P_i\le q_k|Z=z\right\}}{q_k} \times \\
& & \hspace*{68pt} \mathbb{P}_{\vt, \eta}\left\{ C_k^{(i)}|Z=z\right\}d F_Z(z),
\end{eqnarray*}
where the random events $A_k^{(i)}$ and $C_k^{(i)}$ are defined in the proof of Theorem \ref{thm1}.
Moreover, making use of the notation $C_{\bY;k}^{(i,z)}$ introduced in the proof of Theorem \ref{thm1},
we can express $\FDR_{\vt, \eta}(\vp^{LSU})$ by
\begin{eqnarray*}
\FDR_{\vt, \eta}(\vp^{LSU})&=&\sum_{i=1}^{m_0} \frac{q}{m} \sum_{k=1}^{m}
\int \frac{\mathbb{P}_{\vt, \eta}\left\{P_i\le \exp\left(-z \psi^{-1}\left(F_{P_i}\left(\frac{kq}{m}\right)\right)\right)\right\}}{q_k} \times\\
& & \hspace*{68pt} \mathbb{P}_{\vt, \eta}\left\{ C_{\bY;k}^{(i,z)}\right\}d F_Z(z)\\
&=&\sum_{i=1}^{m_0} \frac{q}{m} \sum_{k=1}^{m}
\int \frac{\exp\left(-z \psi^{-1}\left(q_k\right)\right)}{q_k} \mathbb{P}_{\vt, \eta}\left\{ C_{\bY;k}^{(i,z)}\right\}d F_Z(z)\\
&\ge& \sum_{i=1}^{m_0} \frac{q}{m} \sum_{k=1}^{m}
\int \min_{k \in \{1,\hdots,m\}} \left\{\frac{\exp\left(-z \psi^{-1}\left(q_k\right)\right)}{q_k}\right\} \times\\
& & \hspace*{68pt} \mathbb{P}_{\vt, \eta}\left\{ C_{\bY;k}^{(i,z)}\right\}d F_Z(z),
\end{eqnarray*}
where the latter inequality follows from $Y_\ell \sim \text{UNI}[0,1]$ for all $1 \leq \ell \leq m$ and the fact that each $H_{i}$ is a true null hypothesis.

Now, it holds that
\begin{eqnarray*}
\FDR_{\vt, \eta}(\vp^{LSU})&\ge &\sum_{i=1}^{m_0} \frac{q}{m} \int \min_{k \in \{1,\hdots,m\}} \left\{\frac{\exp\left(-z\psi^{-1}\left(\frac{kq}{m}\right)\right)}{kq/m}\right\} \times\\
& & \hspace*{64pt} \sum_{k=1}^{m}  \mathbb{P}_{\vt, \eta}\left\{C_{\bY;k}^{(i,z)}\right\}dF_{Z}(z)\\
&=&\sum_{i=1}^{m_0} \frac{q}{m} \int \min_{k \in \{1,\hdots,m\}} \left\{\frac{\exp\left(-z\psi^{-1}\left(\frac{kq}{m}\right)\right)}{kq/m}\right\} dF_{Z}(z)\\
&=& \int\min_{k \in \{1,\hdots,m\}} \left\{\frac{\exp\left(-z\psi^{-1}\left(\frac{kq}{m}\right)\right)}{kq/m}\right\}  dF_{Z}(z) \sum_{i=1}^{m_0} \frac{q}{m}\\
&=& \int \min_{k \in \{1,\hdots,m\}} \left\{\frac{\exp\left(-z\psi^{-1}\left(\frac{kq}{m}\right)\right)}{kq/m}\right\}  dF_{Z}(z) \frac{m_0 q}{m}.
\end{eqnarray*}

This completes the proof of the theorem.\qed

\subsection*{Proof of Theorem \ref{thm-ex}}
 We plug \eqref{equation-21} into \eqref{th1_eq1} and obtain
\begin{eqnarray*}
F_{\tilde{P}_1,\hdots,\tilde{P}_{m}}(p_1,\hdots,p_m)&=&
\int \prod_{i=1}^m \exp\left(-z \psi^{-1}(p_i)\right) dF_Z(z)\\
&=& \int \exp\left(-z\sum_{i=1}^m \psi^{-1}\left(F_{\tilde{P}_i}(p_i)\right)\right) dF_Z(z)\\
&=& \psi\left(\sum_{i=1}^m \psi^{-1}\left(F_{\tilde{P}_i}(p_i)\right)\right),
\end{eqnarray*}
since the last integral is the Laplace transform of $Z$ at $\sum_{i=1}^m \psi^{-1}\left(F_{\tilde{P}_i}(p_i)\right)$.
Noticing that $P_1,\hdots,P_{m}$ are obtained by componentwise increasing transformations of $\tilde{P}_1,\hdots,\tilde{P}_m$ we conclude the assertion. \qed

\section*{Acknowledgments}
This research is partly supported by the Deutsche Forschungsgemeinschaft via the Research Unit FOR 1735 "Structural Inference in Statistics: Adaptation and Efficiency". 

\bibliographystyle{imsart-nameyear}
\bibliography{FDR-AC}
\end{document}